\documentclass[11pt,oneside]{amsart}

\usepackage[a4paper, total={6in, 8.8in}]{geometry}

\usepackage{amsmath, amsthm, amssymb}    
\usepackage{graphicx,color}                    
\usepackage{hyperref}   
\usepackage{bm}

\newcommand{\bbar}{\begin{pmatrix}}
\newcommand{\ebar}{\end{pmatrix}}
\newcommand{\bmat}{\left[ \begin{array}}
\newcommand{\emat}{\end{array}\right]}
\newcommand{\beq}{\begin{equation}}
\newcommand{\eeq}{\end{equation}}

\newcommand{\bsigma}{{\bm \sigma}}

\newcommand{\Res}{\operatorname{Res}}
\newcommand{\Err}{\operatorname{Err}}

\newtheorem{proposition}{Proposition}
\newtheorem{definition}{Definition}

\title{Approximation and interactive design with exact 3D elastic curves}

\address{Department of Applied Mathematics and Computer Science, Technical University of Denmark \\
Richard Petersens Plads, Building 324,  DK-2800, Kgs. Lyngby,  Denmark}
\author{David Brander}
\email{dbra@dtu.dk}
\author{Jens Gravesen}
\email{jgra@dtu.dk}
\author{Marc Isern}
\email{misha@dtu.dk}

\begin{document}

\begin{abstract}
    An elastic space curve is a critical point of the bending energy subject to appropriate constraints. An analytic representation, equivalent to the spherical pendulum equation, leads to an 11-parameter description of the space of 3D elastic curve segments. We give a numerically stable method for recovering the 11 parameters from a given elastic curve segment. Using this, we give a fast and stable method to approximate an arbitrary space curve segment by a 3D elastica. Applications include interactive design with exact elastic curves and CAD surface rationalization for robotic hot-blade cutting.
\end{abstract}
% Keywords
% MSC2020 primary and secondary classifications

\maketitle

%\begin{figure*}
 %\includegraphics[height=0.19\textheight]{images/nex1c}     \quad \quad
 % \includegraphics[height=0.2\textheight]{images/dini4.png}  
 %\centering
  %\caption{Examples of elastica swept surfaces.}
 %\label{fig:teaser}
%\end{figure*}

\section{Introduction}
Elastic curves, or elastica, are critical points of the bending energy  $\int \kappa^2 \mathrm{d}s$ under length and boundary constraints. They arise naturally when a curve is produced by a flexible physical element, such as a rod, strip, wire or blade. This leads to  many practical applications, such as in cable modeling and routing \cite{jin2022robotic,monguzzi2024potential},  surgical sutures \cite{wang2017real},
curve-completion models in computer vision and image inpainting \cite{mumford1994elastica,shen2003euler},  robotic manipulation  of flexible linear objects  \cite{moll2006path,levin2025dual}, robotic hot-blade cutting \cite{sondergaard2016robotic}, and in geometric design, for example in elastica-ruled surfaces and kinetic structures assembled from elastica strips \cite{lee2020ruled,bi2023design}.

 The basic ingredient needed for many of these problems is to compute an elastic curve that matches a given target shape, in some cases only roughly. Existing methods for computing an elastic curve involve some version of solving a boundary-value problem, such as prescribing end-points, end-tangents and length, and using an optimization to minimize the elastic energy (see Section \ref{related}). When such constraints are paramount this makes sense, but such methods are less well suited to the inverse fitting that we wish to consider, as to adapt this to the problem of matching a \emph{target shape} involves solving many smaller boundary-value problems to reach an interpolation or approximation.  While this can be done even for a large family of curves, it is time-consuming and the outcome (depending on the approach) may be curves with low elastic energy, but not necessarily true elastica.
Our goal is to develop a fast method for shape approximation using exact 3D elastic curves.

In this article we restrict ourselves to twist-free elastic curves.   With this assumption, the shape of the curve is material-agnostic as it is independent of the material-specific torsional response, allowing for a broad range of applications.  An extension of our method is expected to be possible for elastic rods with torsion (see Section \ref{conclusion} below).

\subsection{Approach and Main Results} 
The problem of approximating a \emph{plane} curve by an elastic curve was studied in 
\cite{bgn2016}. The method was based on a parameterization of the space of planar elastica
in terms of elliptic functions. The method we will use for 3D elastica follows a similar logic, the main difference being that the analytic description of the solution space is considerably more complicated.

We represent the \textbf{space of 3D elastic curve segments} as an 11-dimensional parameter space; that is, any elastic curve segment is (essentially uniquely) determined by an 11-tuple of numbers, $\bsigma = (\sigma_1, \dots \sigma_{11})$.  (Only two of these parameters are intrinsic elastica parameters; the remaining parameters describe scale, translation, rotation, and the choice of segment on the normalized model curve. See Section \ref{parameters}).  
\begin{figure}[htb]   
\centering
 \begin{tabular}{cc}
   \includegraphics[height=0.18\textheight]{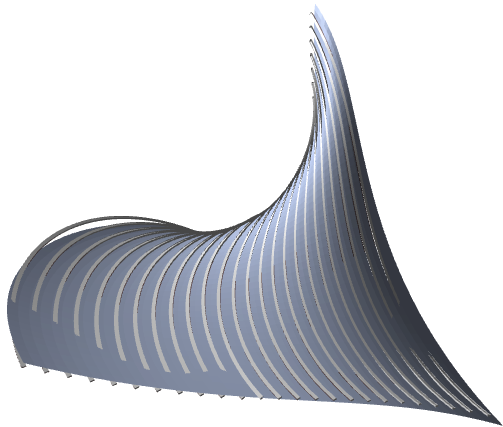} \quad
& \quad \includegraphics[height=0.18\textheight]{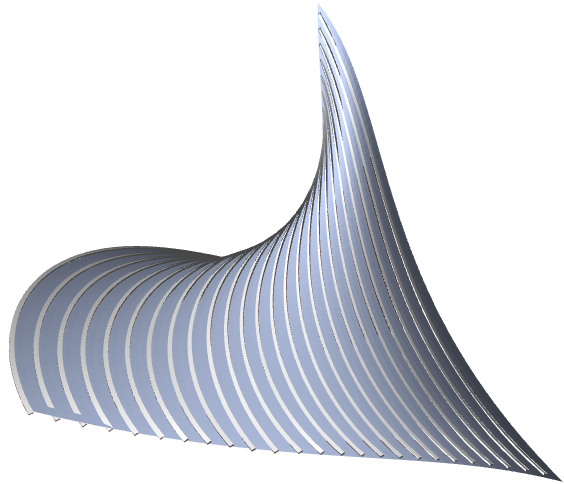}
\end{tabular}
 %%%%%%%%%%%%%%%%%%%%%%%%%
 \caption{A surface (blue) approximated by a family of elastic space curves (grey strips). 
 Left: fast least-squares method, maximum error $0.013$. Right: with  subsequent nonlinear optimization, maximum error $0.0038$.}
    \label{FigSurfApprox}
\end{figure}

\textbf{The least-squares elastica $\Gamma(t;\bsigma_{\gamma})$: } A space curve
 $\gamma(t)$ is an elastic curve if and only if it satisfies a 3rd order nonlinear equation, involving the relevant elastic curve parameters in $\bsigma$. 
Solving this equation for the \emph{parameters} (rather than the solution itself),  and some other relevant geometric equations of the curve, in the least-squares sense leads to a series of linear least-square problems. We obtain an $11$-tuple $\bsigma_\gamma$, the parameters for an elastic curve, which we designate $\Gamma(t;\bsigma_{\gamma})$.  If $\gamma$ was itself an elastic curve, then $\bsigma_\gamma$ is the set of parameters for $\gamma$. % Otherwise, if $\gamma$ is a twice differentiable curve, then $\Gamma(t;\bsigma_{\gamma})$ is close in the $L^2$ sense to $\gamma$. %%this is wrong
Moreover, the first residual $\Res_\lambda$ in
 the sequence of least-squares problems correlates well with the $L^2$ distance between $\gamma$ and its fitted elastica $\Gamma(t;\bsigma_{\gamma})$.   

\textbf{Application 1: interactive  modeling with exact 3D elastica and elastic splines.} Since $\bsigma_\gamma$ and $\Res_\lambda$ can be computed very quickly by solving only linear problems, and $\Gamma(t;\bsigma_{\gamma})$ is an actual elastic curve, whose deviation from a target shape $\gamma$ can be roughly estimated by $\Res_\lambda$, these are  ingredients
 that can be used for dynamic geometric modeling with 3D elastic curves.  
 This is explained in Sections \ref{Design},
 \ref{Residual} and
 \ref{splinesection}.
 
 \textbf{Application 2: approximation of target CAD curves by exact elastic curves.}
An important motivation for this work is the problem of surface rationalization for robotic hot-blade cutting. This is a generalization of hot-wire cutting of EPS foam,
where, instead of a straight wire, a heated elastic rod is used (see \cite{sondergaard2016robotic}).
In a project linked to this work, we seek to develop methods for surface rationalization using foliations by continuously varying 3D elastic curves. The basic problem here is to approximate a target curve as closely as possible by an elastic curve.  To solve this problem, we take the least-squares elastica $\Gamma(s;\bsigma_{\gamma})$ and then apply a gradient-driven optimization to move inside the 11 parameter space to an optimized solution. This is discussed in Section \ref{Approximation}. An example of the result is shown in Figure \ref{FigSurfApprox}.  

\subsection{Structure of this paper}
In Section \ref{related} we give a brief survey of related work.
In Sections \ref{var} and \ref{parameters} we describe the analytic model of 3D elastic curves, including a derivation of the relevant equations based on the spherical pendulum equation, the representation of all solutions in terms of elliptic functions, and the 11-dimensional parameter space of elastic curves.

The least-squares method for recovering the control parameters for an elastic curve is given in Section \ref{lssection}. This is the most important result, as it can be applied to any curve to obtain an elastic curve with a similar shape.  As already mentioned this method can be utilized for various problems that require exact elastic curves guided by the global shape of a target curve, and that is described in the remainder of the paper.

\section{Some Background and Related Work}\label{related}
We are not aware of any work directly addressing the problem of finding a 3D elastic curve that is close in shape to a given space curve segment.   There are, however, many works on what may be termed the \emph{forward} problem for elastic curves, namely to compute the elastic curve that satisfies some constraints, and we will try to summarize this here.

\subsection{Elastic Curves}
Euler's elastica appeared in design before the introduction of computers, because these curves are the natural shape for the wooden splines that were used in drafting.
 Elastic curves have long been used in minimum-energy interpolation \cite{horn1983curve}, while more recent active-bending work studies how planar elastic strips of variable bending stiffness can be designed so that their equilibrium shapes realize planar target curves \cite{hafner2021design}.

In the planar Euler--Bernoulli setting, the  elastica family can  be written explicitly in terms of Jacobi elliptic functions, and modern expositions show how solutions can be constructed under several endpoint conditions \cite{arroyo2020boundary}. These methods are relevant here mainly as a point of contrast: even when the forward family is explicit, solving the associated boundary-value problem still requires branch selection and nontrivial parameter recovery, whereas our goal is not to satisfy prescribed boundary conditions but to fit an exact elastic curve directly to a target shape.

\subsection{Rod theories}
Rod theory and beam models allow for more deformation modes than the twist-free elastic-curve setting considered here. Reviews of rod-based and soft-robot modeling can be found in \cite{alessi2024rod,armanini2023soft}. They distinguish Euler--Bernoulli and Kirchhoff descriptions from Timoshenko and Cosserat theories, where additional generality is obtained by introducing shear, extension, twist, and full director dynamics. Elastic curves appear as the inextensible twist-free rods with no shear. Representative examples include spline-based Kirchhoff formulations for simulation of sutures and cables \cite{panneerselvam2020constrained}, as well as discrete differential-geometric and Cosserat-style frameworks such as {\fontfamily{cmss}\selectfont DisMech} and {\fontfamily{cmss}\selectfont Elastica} \cite{choi2024dismech,armanini2023soft}. These methods are useful for contact or dynamic problems, but they address a different computational regime from the inverse fitting problem studied here.

From a numerical perspective, rod and elastica problems are commonly approached in several different ways. Elastica can be written in closed form using Jacobi elliptic functions and elliptic integrals, and so the main computational task becomes numerical parameter recovery and evaluation of these special functions, which can be done fast (see \cite{carlson1995numerical} and for an implementation \cite{batista2019elfun18}).  Other approaches solve the governing Kirchhoff or elastica equations directly by ODE or boundary-value integration; representative comparisons of ODE-based Kirchhoff integration, Euler elastica boundary-value formulations, and discrete elastic rods appear for example in \cite{tola2022simulation}. Discretization-based methods include geometrically exact finite-element formulations for Bernoulli--Euler or Kirchhoff rods \cite{da2020simple}, as well as discrete differential-geometric models such as Discrete Elastic Rods (DER) \cite{bergou2008discrete}. These methods are highly effective when the goal is to simulate rod mechanics under specified loads and constraints. By contrast, our setting is an inverse geometric fitting problem: rather than numerically solving a rod model for each query, we recover nearby \emph{exact} elastic curves.

\subsection{Bending-active form finding}
A related but still broader strand appears in the computational form finding of bending-active structures. Comparative studies of dynamic-relaxation-based tools and related beam models show that both accuracy and runtime depend strongly on the numerical scheme, discretization, and underlying beam theory \cite{cuvilliers2018comparison,lazaro2018mechanical}. We cite this literature only to contextualize computational trade-offs: it addresses structural simulation and equilibrium finding rather than fast inverse fitting of exact elastic-curve segments to a given input curve.\\

%%%%%%%%%%%%%%%%%%%%%%%%%%%%%%%%%%%%%%%%%%%%%%%%%%%%%%%%%%%
\section{Variational Approach and Analytic Representation}\label{var}
The theory of elastic curves is classical and good accounts of it can be found in \cite{singer2008lectures} and \cite{pinkall2024differential}. In this section we describe the variational characterization of elastic curves using a different approach, based on the fact that the Euler-Lagrange equations of an elastic curve coincide with the equations of motion of a spherical pendulum.  We conclude this section by deriving the analytic representation in terms of elliptic functions. Equivalent formulas can be found, for instance, in \cite{ameline2017classifications}.

\subsection{Variational formulation} \hfill\\
An \textit{elastic curve} $\gamma$ is a critical point of the bending energy functional $$\mathcal{B}(\gamma)=\frac{1}{2}\int_{\gamma}\kappa^2 \mathrm{d}s$$ subject to length and boundary constraints. In addition, we choose the curve to be unit-speed parameterized.  %The image of an elastica in $\mathbb{R}^n$ lies in a three-dimensional subspace; so it suffices to look at $n=3$ to classify all elastica (cf. \cite{miura2024elastic}). \\

 For a (sufficiently smooth) unit-speed parameterized curve $\gamma : [0,L] \rightarrow \mathbb{R}^3 $ this variational problem amounts to finding critical points of
\begin{equation*}
    \mathcal{E}(\gamma)=\mathcal{B}(\gamma)+\Big\langle\lambda,\gamma(0)+\int_0^L\gamma'(s)\mathrm{d}s-\gamma(L) \Big\rangle
\end{equation*} 
for some Lagrange multiplier $\lambda\in \mathbb{R}^3$, which encodes the constraint of fixed endpoints.

\subsection{Euler-Lagrange equations and the spherical pendulum}\hfill\\
As the curve is unit-speed parameterized, its derivative $\gamma'$ takes values on the unit sphere. Thus, a natural approach is to use spherical coordinates for the tangent vector; using the azimuthal angle $\theta$ and the polar angle $\phi$ we have 
\[\gamma'(s)=(\cos\theta(s)\sin\phi(s),\sin\theta(s)\sin\phi(s), \cos\phi(s)).\]

We consider a variation $\gamma_t$ of $\gamma$ by considering variations of the angle functions. For each $t$, $\gamma_t$ is the curve with the same startpoint as $\gamma$ and with angle functions $\theta_t=\theta+t\varphi$ and $\phi_t=\phi+t\eta$, for differentiable functions $\varphi, \eta$ with $\varphi(0)=\varphi(L)=\eta(0)=\eta(L)=0$. Then the end tangents of $\gamma_t$ agree with those of $\gamma=\gamma_0$ and the curves have the same length.\\

The first variation of $ \mathcal{E}$ can be computed in the standard way using integration by parts (using $\varphi(0)=\varphi(L)=\eta(0)=\eta(L)=0$ to kill the boundary terms) as:

\begin{align*}
     \frac{\text{d}\mathcal{E}(\gamma_t)}{\text{d}t}\Big|_{t=0}
    =&\int_0^L\eta\left(\sin\phi\cos\phi \;(\theta')^2-\phi''+\lambda_{\theta}\right)\text{d}s \\
    &- \int_0^L\varphi\sin\phi\left(2\theta'\phi'\cos\phi+\sin\phi\;\theta''- \lambda_{\phi}\right)\text{d}s,
\end{align*}
where $\lambda_\theta=-\lambda_1\sin\theta+\lambda_2\cos\theta$ and $\lambda_\phi=\lambda_1\cos\theta\cos\phi+\lambda_2\sin\theta\cos\phi-\lambda_3\sin\phi$ are the spherical components of $\lambda= \lambda_r\hat{\underline{\mathbf{r}}}+\lambda_\theta \hat{\underline{\mathbf{\theta}}}+\lambda_\phi \hat{\underline{\mathbf{\phi}}}$.\\

Since $t=0$ is a stationary point of $t\mapsto\mathcal{E}(\gamma_t)$ and $\varphi, \eta$ are arbitrary, we conclude that the angle functions of an elastic curve must satisfy:
\begin{align*}
    \phi''-\sin\phi\cos\phi \;(\theta')^2 =\lambda_\phi,\\
    2\theta'\phi'\cos\phi+\sin\phi\;\theta''=\lambda_\theta.
\end{align*}
These are precisely the equations for a spherical pendulum with gravity vector $\lambda$. That is, the tangent vector $\gamma'$, which takes values on the unit sphere, evolves as a spherical pendulum subject to a constant gravity vector $\lambda$.\\

Since  $\gamma'''= -\kappa^2\hat{\underline{\mathbf{r}}}+(2\theta'\phi'\cos\phi\sin\phi+\theta''\sin^2\phi) \hat{\underline{\mathbf{\theta}}}+(\phi''-(\theta')^2\sin\phi\cos\phi) \hat{\underline{\mathbf{\phi}}}$,
 we can rewrite the equation of the pendulum independently of our choice of coordinates as is stated in the following proposition.
\begin{proposition}[cf.\cite{pinkall2024differential}]  \label{elasticprop}
    Let $\gamma:[0,L]\rightarrow\mathbb{R}^3$ be a unit-speed parameterized curve. Then $\gamma$ is an elastic curve iff 
    \begin{align}
         \gamma'''- \langle\gamma''', \gamma'\rangle\gamma'=\lambda- \langle\lambda, \gamma'\rangle\gamma' ,\label{torque}
    \end{align}
    for some constant $\lambda\in\mathbb{R}^3$. Equivalently,
    \begin{align}
    \gamma'''\times \gamma'=\lambda\times\gamma'\,.   \label{physical}
    \end{align}
\end{proposition}
\medskip
 Note that one can also arrive at equation (\ref{physical}) from physical considerations by looking at forces acting on a rod at equilibrium.
 
\subsection{Conserved quantities}\hfill\\
The conserved quantities of a spherical pendulum - its energy and angular momentum about the direction of gravity - will play a central role, as they determine the motion of the pendulum. They are two intrinsic parameters associated to an elastic curve:

\begin{proposition}
    Let $\gamma$ be an elastic curve, and $\lambda$ the 
    constant in Proposition \ref{elasticprop}. Then 
    \begin{align*}
        \lambda_0= \frac{1}{2} \langle \gamma'', \gamma'' \rangle - \langle \lambda, \gamma' \rangle  \quad \hbox{and} \quad \omega=\langle\lambda,\gamma'\times\gamma'' \rangle
    \end{align*}
    are constant along the curve.
\end{proposition}\medskip
The constant $\lambda_0$ is commonly referred to as the $\textit{tension}$ of an elastic curve and corresponds to the energy of the spherical pendulum.
 Using $\lambda_0$ we can rewrite (\ref{torque}) as:
\begin{align}
    \gamma'''+\frac{3}{2}\langle\gamma'',\gamma''\rangle\gamma'-\lambda_0\gamma'=\lambda \, . \label{elastica}
\end{align}

 By definition, the angular momentum of the pendulum is $\Omega=\gamma'\times\gamma''$ and its angular momentum about the direction $\lambda$ is $\omega=\langle\lambda,\Omega\rangle=\langle\lambda,\gamma'\times\gamma'' \rangle$. In particular, $\omega=\langle\lambda,\gamma'\times\gamma'' \rangle=\langle\gamma''',\gamma'\times\gamma''\rangle=\kappa^2\tau$. So remarkably, $\kappa^2\tau$ is a constant in the case of elastic curves. \\
 
 If $\lambda=0$ then $\lambda_0=\frac{1}{2}\kappa^2$; thus, $\kappa$ must be a constant, and $\omega=0$ so $\tau=0$, so either $\gamma$ is a straight line or a circle. In the following, assume $\lambda \neq 0$.\\

Equation \eqref{physical} can be integrated to $\gamma''\times \gamma'=\lambda\times\gamma+A$, for some constant $A\in\mathbb{R}^3$. We may choose the origin of our coordinate system such that $A$ is parallel to $\lambda$. And since $\langle\lambda, A\rangle=\langle\lambda, \gamma''\times\gamma'\rangle=-\omega$ we can write
\begin{equation}
    \gamma''\times \gamma'=\lambda\times\gamma-\frac{\omega}{||\lambda||^2}\lambda \,. \label{physics}
\end{equation}

\subsection{Normalized elastica}\hfill\\
When $\lambda\neq0$, up to appropriate rotation and scaling, we may choose $\lambda=(0,0,-1)$. In that case $\lambda_0=\frac{1}{2}((\theta'\sin\phi)^2+(\phi')^2)+\cos\phi$ and $\omega=-\theta'\sin^2\phi$. \\ Let $z(s)=\gamma_3(s)$ denote the vertical coordinate. Then  $\cos\phi=z'$ and combining both conserved quantities one arrives at a differential equation for $z'$, which can be solved using elliptic functions:
%We will call elastica with $\lambda=(0,0,1)$ normalized elastica. A normalized elastica whose tangent satisfies the pendulum equation with energy $\lambda_0$ and angular momentum about gravity $\omega$ will be denoted by $\zeta_{\lambda_0,\omega}$.We will give an explicit parametrization of normalized elastica in terms of elliptic functions. To do this it is convenient to use cylindrical coordinates $(r,\Theta,z)$. Then $\cos\phi=z'$, with which one obtains the following equation:
\begin{align}
    (z'')^2=2(\lambda_0-z')(1-(z')^2)-\omega^2.
    \label{poly}
\end{align}
Let $\beta_1\geq1\geq\beta_2\geq\beta_3\geq-1$ be the three roots of the cubic polynomial $$P(x)=2(\lambda_0-x)(1-x^2)-\omega^2,$$ and define  
\begin{align*}
    m=\frac{\beta_2-\beta_3}{\beta_1-\beta_3}, && \delta=\sqrt{\frac{\beta_1-\beta_3}{2}}.
\end{align*}
Then, we have
    \begin{align}
    z'= \beta_3+(\beta_2-\beta_3)\mathrm{sn}^2(\delta s|m), \label{height}
    \end{align}
    where $\mathrm{sn}$ is the elliptic sine function. For the notation of elliptic functions, we use the Jacobi normal form (cf. \cite{batista2019elfun18}).  \\
    
    Note that $z'\in[\beta_3,\beta_2]$, so $\beta_2$ and $\beta_3$ are the maximum and minimum heights of the pendulum. For  (\ref{poly}) to have a solution, there must exist $x\in[-1,1]$ with $P(x)\geq0$. This is the case if the local maximum $x_{\text{max}}$ satisfies $P(x_{\text{max}})\geq 0$. This inequality in terms of the conserved quantities $\lambda_0$ and $\omega$ is
\begin{align}
    \omega^2\leq \frac{4}{27}((\lambda_0^2+3)^{3/2}-(\lambda_0^3-9\lambda_0)), \label{E-L-condition}
\end{align} which defines the admissible domain of pairs $(\lambda_0,\omega)$, together with the inequality that follows from the definition $\lambda_0\geq -1$.\\

\begin{figure}[htb]
    \centering
    \includegraphics[width=0.5\linewidth]{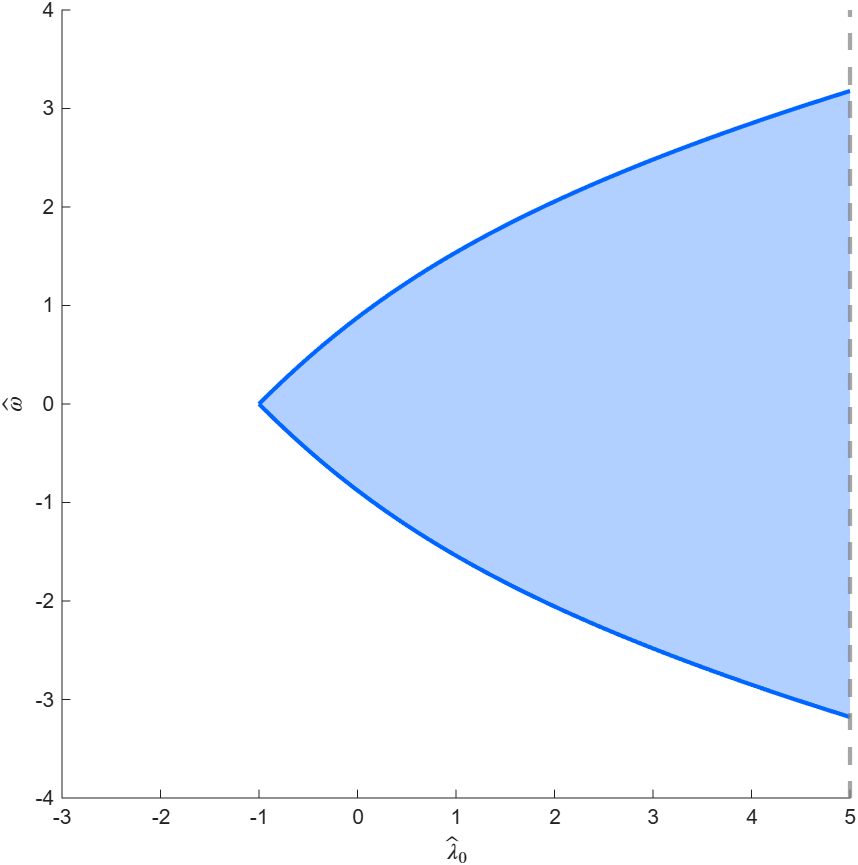}
    \caption{Admissible domain for $(\hat\lambda_0,\hat\omega)$}
    \label{fig:placeholder}
\end{figure}

Equation (\ref{height}) can be integrated to find:
\[
    z= z_0+ \beta_1 s - 2\delta E(\delta s| m),
\]
for some constant $z_0$, where $E$ is the second incomplete elliptic integral.\\

By taking the dot product of $\gamma'$ and (\ref{physics}) on both sides we find:
\begin{align*}
      0 & = \langle\gamma\times \gamma', \lambda\rangle  - \frac{\omega}{||\lambda||^2} \langle\lambda, \gamma'\rangle \,, 
\end{align*}
 and the dot product of $\gamma'$ with the cross product of (\ref{physics}) and $\lambda$ we find:
\begin{align*}
    \langle \gamma'',\lambda\rangle &=  \langle\gamma, \gamma'\rangle- \langle\lambda,\gamma\rangle\langle\lambda, \gamma'\rangle \, .
\end{align*}
In cylindrical coordinates $\gamma(s)=(r(s)\cos\Theta(s),r(s)\sin\Theta(s),z(s))$, for $\lambda=(0,0,-1)$, these equations are:
\begin{align}
    0  & =  -r^2\Theta'+\omega z' \label{cylang},\\
    -z'' & = r'r \label{radius}.
\end{align}
Integrate (\ref{radius}) to find $r^2$. The integration constant can be determined using $\dot{r}^2+r^2\dot{\theta}^2+\dot{z}^2=1$ and (\ref{poly}) to be $2\lambda_0-\omega^2$. Subsequently, we use (\ref{cylang}) to solve for $\Theta$. This yields the following analytic representation of an elastic curve in cylindrical coordinates:
    \begin{align*}
    z&= z_0+ \beta_1 s - 2\delta E(\delta s|m),\\
r^2&=2\lambda_0-\omega^2-2\beta_3-2(\beta_2-\beta_3)\text{sn}^2(\delta s|m), \\
    \Theta&=\Theta_0-\frac{\omega}{2}s+\Pi_{\text{ext}}(\lambda_0,\omega,\delta s),
\end{align*}
where $E$ is the elliptic integral of the second kind. The angle function is not always defined, so this requires a small convention. In the case of a straight line (so $(\lambda_0,\omega)=(-1,0)$) it is never defined. Otherwise, let
 \begin{align*}
     n=\frac{2(\beta_2-\beta_3)}{2\lambda_0-\omega^2-2\beta_3}
 \end{align*} and define
\begin{align}
   \Pi_{\text{ext}}(\lambda_0,\omega,\delta s):=\frac{\omega(2\lambda_0-\omega^2)}{2\delta(2\lambda_0-\omega^2-2\beta_3)}\Pi(n;\delta s| m), \label{extension}
\end{align}
where $\Pi$ is the incomplete elliptic integral of the third kind. The first factor removes the singularities of $\Pi$, so although $\Pi$ may diverge, the formula \eqref{extension} always remains finite.\\
The elliptic integral $\Pi$ may have singularities for special parameter values, in particular when $n=1$, which here corresponds to the radius reaching zero. So the apparent singularity is only a singularity of the cylindrical coordinate representation, not of the elastic curve itself. For fixed $m$, in the limit $n \rightarrow 1$,  $\Pi_{\text{ext}}$ converges to a step function, corresponding to the limiting jump in the cylindrical angle as the curve passes through the axis.

\begin{definition}
    For admissible $(\lambda_0,\omega)$, the corresponding normalized elastica, denoted by $\zeta_{\lambda_0,\omega}$, is defined by the above formulas with $z_0=0$ and $\Theta_0=0$.
\end{definition}
\begin{figure}[ht] 
    \centering
    \includegraphics[width=0.9\linewidth]{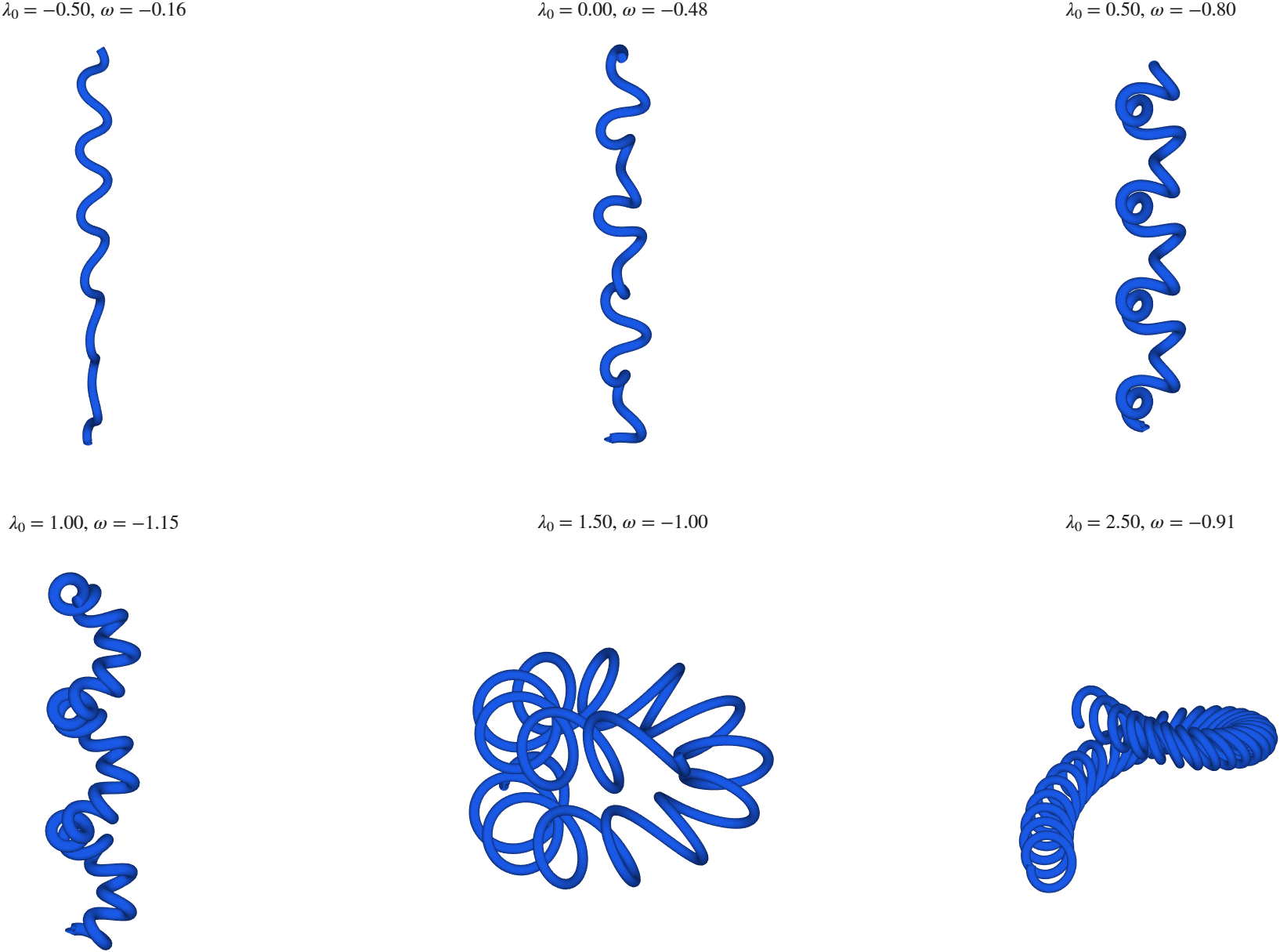}  %
    \caption{Generic normalized elastic curves for different values of $(\lambda_0,\omega)$ for symmetric ranges of $s$ about zero}
    \label{FigExamples}
\end{figure}

The analytic formulas above also give us a nice geometric picture.  Generically, an elastic curve winds around a core helix and is contained in a tube-shaped envelope with circular cross-section, see Figure \ref{FigExamples}. The core helix winds itself around the direction of $\lambda$. The intrinsic parameters of the normalized elastica therefore determine not only the curve, but can be used to also determine the radius, tube width, and pitch of this helical envelope (see \cite{ameline2017classifications}). Spatial elastica also include special cases such as straight lines, circles, planar elastica and helices, arising from degeneracies of the helical envelope.

%%%%%%%%%%%%%%%%%%%%%%%%%%%%%%%%%%%%%%%%%%%%%%%%%%%%%%%%%%%%%%%%%
\section{Parametrization of the space of elastica}  \label{parameters}
As seen in the previous section, every elastic curve (except circles arising from $\lambda=0$) can, up to scaling and a Euclidean motion, be represented by a normalized elastic curve. Thus,  an arbitrary segment of an elastica is determined in two steps. Choose $\hat{\lambda}_0$ and $\hat{\omega}$ to pick the type of normalized elastica $\zeta_{\hat{\lambda}_0, \hat{\omega}}$ and pick a segment by specifying a starting point $s_0$ and an endpoint $s_0+\ell$. Second, apply a scaling  factor $S>0$, and a rotation given by $R\in SO(3)$ and a translation  by $(x_0,y_0,z_0)$.  \\
Note that $\ell$, $\hat{\lambda}_0$ and $\hat{\omega}$ are the length, tension and angular momentum of the normalized elastica. The scaled curve $\gamma$ has length $L=S\ell$, tension $\lambda_0=S^{-2}\hat{\lambda}_0$ and angular momentum $\omega=S^{-3}\hat{\omega}$. The scaling relation for $\lambda$ is $\|\lambda\|=S^{-2}$.
 \\
\begin{figure}[h]
    \centering
\includegraphics[height=0.2\textheight]{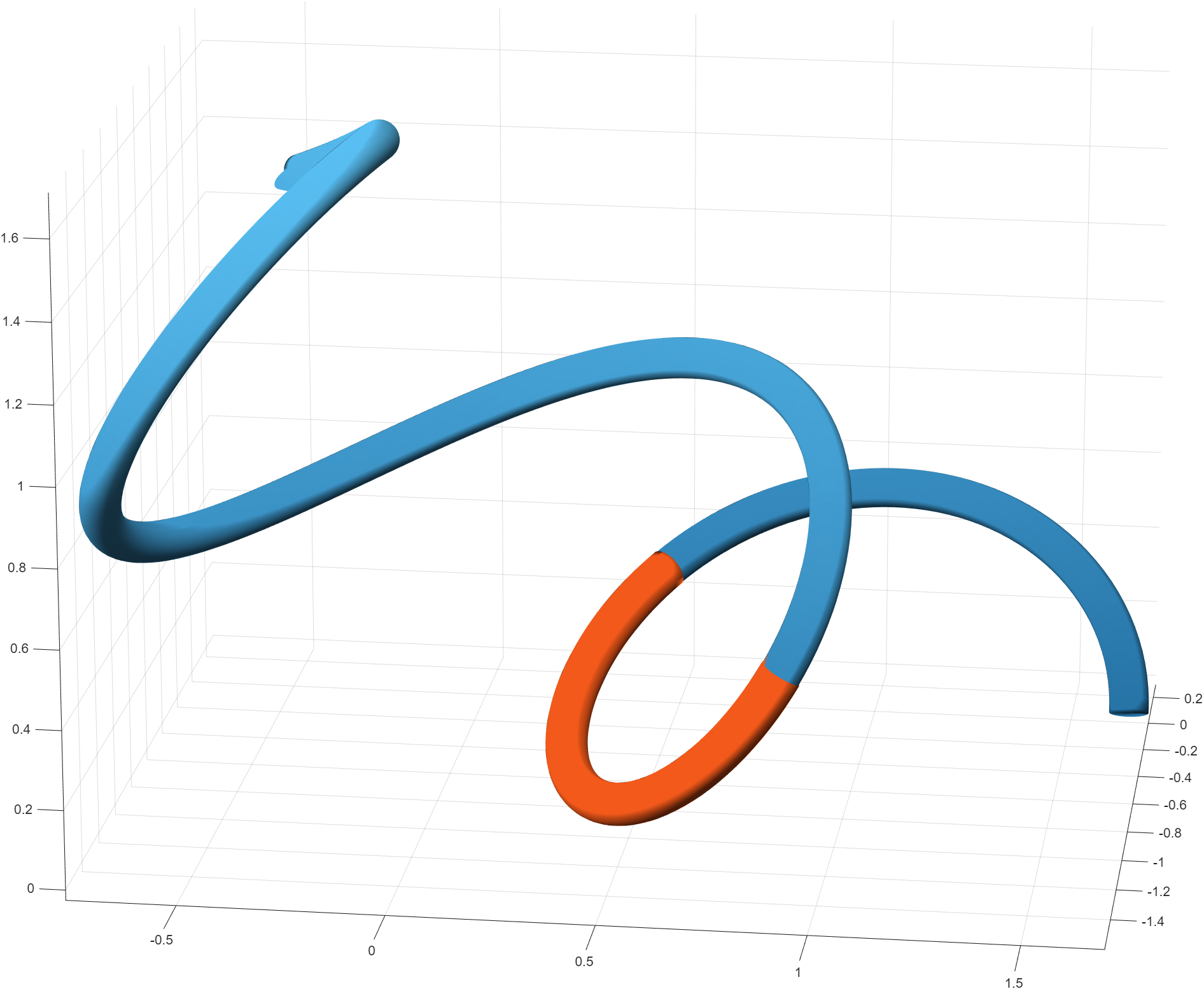}
     \caption{Normalized elastic curve $\zeta_{1,1}$ in blue, with a segment with $s_0=1.6$ and $\ell=1.9$ in red.}
\end{figure}
In summary, we have a unit-speed parameterized curve:
\begin{align}
    \gamma(s)=SR\zeta_{\hat{\lambda}_0,\hat{\omega}}(s/S+s_0) + (x_0,y_0,z_0) \text{  for } s\in[0,L].
\end{align}
It is convenient (for instance when applying an optimization) to make the domain independent of the parameters, so one can reparameterize the curve in the interval $[0,1]$:
\begin{align}
   \Gamma(t,\bsigma)= \gamma(t)=SR\zeta_{\hat{\lambda}_0,\hat{\omega}}(\ell t+s_0) + (x_0,y_0,z_0) \text{  for } t\in[0,1].
\end{align}

Thus, elastic curve segments are described by eleven parameters:
\[
 \bsigma=( \hat{\lambda}_0, \hat{\omega}, s_0, \ell,
  S, R, (x_0,y_0,z_0))
\]

Note that $\operatorname{sign}(\omega)$ determines the chirality (or planarity when $\omega=0$) of the curve. We may restrict to $\omega\geq 0$ by allowing reflections, so $R\in O(3)$.\\
The parameters $\hat\lambda_0$ and $\hat\omega$ are nice because they have a physical meaning via the pendulum equation satisfied by $\gamma'$. However, the admissibility domain is not compact. 
An alternative formulation is to use the equivalent intrinsic parameters given in \cite{singer2008lectures}. Instead of $\hat\lambda_0$, $\hat\omega$ and $S$, they use $\kappa_0=S^{-1}\sqrt{\hat\lambda_0-\beta_1}$, $w=\sqrt{\frac{\hat\lambda_0-\beta_1}{\beta_1-\beta_3}}$ and the elliptic modulus $k=\sqrt{m}=\sqrt{\frac{\beta_2-\beta_3}{\beta_1-\beta_3}}$. The admissible set for $w$ and $k$ is given by $0\leq k\leq w\leq 1$, which is compact.

%%%%%%%%%%%%%%%%%%%%%%%%%%%%%%%%%%%%%%%%%%%%%%%%%%%%%%%%%%%%%%%%%%%%%%%%%
\section{Finding the control parameters by a least-squares method}  \label{lssection}
Let $\gamma:[0,L]\rightarrow \mathbb{R}^3$ be a sufficiently smooth unit-speed parameterized curve. If $\gamma$ is an elastic curve, its derivatives satisfy the identities derived in Section \ref{var}. For a general input curve, these identities do not hold exactly. We therefore recover the elastic parameters by enforcing these identities in a least-squares sense. This leads to a sequence of parameter estimation steps described in this section. The parameters obtained yield a first guess for the nonlinear optimization described in Section \ref{Approximation}.

Note that if the input curve is, for example, a discrete curve, or a curve of lower differentiability, then we can still apply the method
in this section: we simply first replace the curve with a suitably smooth 
approximation, such as a cubic spline approximation.

\subsection{Recovery of $\lambda_0$ and $\lambda$}\hfill\\
An elastic curve satisfies equation (\ref{elastica}):
\begin{align*}
    \gamma'''+\frac{3}{2}\langle\gamma'',\gamma''\rangle\gamma'-\lambda_0\gamma'=\lambda.
\end{align*}

For a general curve $\gamma=(\gamma_1,\gamma_2,\gamma_3)$ of length $L$ we solve (\ref{elastica}) in the least-squares sense:
\begin{equation*}
    \underset{\lambda_0,\lambda_1,\lambda_2,\lambda_3}{\text{minimize}}\int_{0}^{L}||\gamma'''+\frac{3}{2}\kappa^2 \gamma' - \lambda_0\gamma'-\lambda||^2\,\mathrm{d}s.
\end{equation*}
This leads to the following linear system, writing $\Delta\gamma= (\Delta\gamma_1,\Delta\gamma_2,\Delta\gamma_3):=\gamma(L)-\gamma(0)$:
\begin{equation}
    \begin{pmatrix}
         
        L & \Delta\gamma_1 & \Delta\gamma_2 & \Delta\gamma_3 \\
        \Delta\gamma_1 & L & 0 & 0   \\
        \Delta\gamma_2 &0 & L &  0 \\
        \Delta\gamma_3 & 0 & 0 & L \\
    \end{pmatrix}
     \begin{pmatrix}
        \lambda_0\\
        \lambda_1\\
        \lambda_2 \\
        \lambda_3\\
       
    \end{pmatrix}
    =
     \begin{pmatrix}
        \mathcal{B}(\gamma) \\
        \int_0^{L}(\frac{3}{2}\kappa^2\gamma_1'+\gamma_1''')\,\mathrm{d}s\\
        \int_0^{L}(\frac{3}{2}\kappa^2\gamma_2'+\gamma_2''')\,\mathrm{d}s \\
        \int_0^{L}(\frac{3}{2}\kappa^2\gamma_3'+\gamma_3''')\,\mathrm{d}s \\
    \end{pmatrix}\, , \label{normal equation}
\end{equation}
where we used $||\gamma'||^2=\gamma_1'^2+\gamma_2'^2+\gamma_3'^2=1$ to simplify.\\
A scale-invariant residual for this is :
\begin{align}
  \label{Eres}
\operatorname{Res}_\lambda=   L^3\int_{0}^{L}||\gamma'''+\frac{3}{2}\kappa^2 \gamma' - \lambda_0\gamma'-\lambda||^2\,\mathrm{d}s\,.
\end{align}
Note that the matrix in the normal equation is singular if and only if the curve is a straight line. The normal equations depend on derivatives up to third order. Thus, with respect
to an $H^3$-metric, the recovered parameters are stable under small perturbations of the input curve. More
precisely, whenever the normal equation \eqref{normal equation} is nonsingular, it depends continuously on
    $\gamma^\prime, \gamma^{\prime \prime}$ and $\gamma'''$.

\subsection{Scaling}\hfill\\
Once $\lambda$ has been determined, the scale of the corresponding normalized
elastica is fixed by the relation $\|\lambda\|=S^{-2}$. Thus,
\[
  S=\|\lambda\|^{-1/2}.
\]
The normalized segment length is therefore $\ell=\frac{L}{S}=L\sqrt{\|\lambda\|}$, and the normalized parameters are
$\hat{\lambda}_0=S^2\lambda_0=\frac{\lambda_0}{\|\lambda\|}$ and $  \hat{\omega}=S^3\omega = \frac{\omega}{\|\lambda\|^{3/2}}$.
Hence, the normalized elastica associated with the input curve $\gamma$ is $\zeta_{\hat{\lambda}_0,\hat{\omega}}$, where $\omega$ still has to be recovered.

\subsection{Recovery of $\omega$}\hfill\\
For elastic curves, $\omega=\kappa^2\tau$ is constant, so we get the least-squares solution by averaging:
\begin{align*}
    \bar\omega=\frac{1}{L}\int_0^{L} \kappa^2\tau\,\mathrm{d}s ,
\end{align*}
with scale-invariant residual:
\begin{align*}
    L^5\int_0^L(\kappa^2\tau-\bar\omega)^2\,\mathrm{d}s .
\end{align*}
 Since $\kappa^2\tau=\operatorname{det}(\gamma',\gamma'',\gamma''')$, this is also stable under small perturbations with respect to an $H^3$-metric.\\
The corresponding normalized parameter $S^3\bar{\omega}$ need not be admissible for a general
input curve. We therefore project it onto the admissible interval given by (\ref{E-L-condition}). Thus, we set
\[
  {\omega} =
  \operatorname{sign}(\bar\omega)
  \min\left\{
  |\bar\omega|,S^{-3}\left(\frac{4}{27}\left((\hat{\lambda}_0^2+3)^{3/2}-(\hat{\lambda}_0^3-9\hat{\lambda}_0)\right)\right)^{1/2}\right\}.
\]  
and $\hat\omega=S^3\omega$.

\subsection{Recovery of $s_0$}\hfill\\
To determine a choice of $s_0$ recall the formula for the derivative of the 
$z$-coordinate of
$\zeta_{\hat{\lambda}_0,\hat{\omega}}$ appearing in Equation \eqref{height}, namely:
\begin{align*}
   z'(s)= \beta_3+(\beta_2-\beta_3)\text{sn}^2(\delta (s+s_0),m) , 
\end{align*}
where $\beta_2, \beta_3, \delta$ and $m$ are constants determined from $\hat\lambda_0$ and $\hat\omega$. Note that, $s_0$ is only determined up to adding the period $2K(m)/\delta$ of the elliptic functions.\\
For the target curve $\gamma$ the value corresponding to the derivative of the $z$-coordinate of the normalized elastica at $s$ is $\langle\frac{-\lambda}{||\lambda||},\gamma'(Ss)\rangle$, which following our pendulum analogy is the height of the pendulum under a gravity vector $\lambda$.\\
Then $s_0$ is determined by the one-dimensional periodic least-squares problem:
\begin{align}
    s_0=\underset{a\in \left[0, \frac{2K(m)}{\delta}\right)}{\operatorname{argmin}}\int_0^\ell\Big(\Big\langle\frac{-\lambda}{||\lambda||},\gamma'(Ss)\Big\rangle-\beta_3-(\beta_2-\beta_3)\text{sn}^2(\delta(s+a)|m)\Big)^2\mathrm{d}s. \label{stables_0}
\end{align}
(when $m=1$, $K(m)$ diverges and the period becomes infinite, in this limiting case the minimization is taken over $\mathbb{R}$). Away from the degenerate case when $\beta_2=\beta_3$ (i.e. $m=0$), where the minimizer is not isolated, this least-squares solution depends on the first derivative $\gamma'$ and on parameters recovered in a stable manner under $H^3$ noise, so it is itself also stable under $H^3$ noise.\\

The periodic problem \ref{stables_0} gives a stable way to recover $s_0$, however, \ref{stables_0} is not linear. If a purely linear least-squares step is desired, one may do so by using the inverse of \ref{height}:
\begin{align}
    s_0=\frac{1}{\delta}  F\left(\sqrt{\frac{ z'-\beta_3}{\beta_2-\beta_3}}\Big|m\right)-s , \label{inverse}
\end{align}
where $F$ is the elliptic integral of the first kind. This is only defined when $\beta_2\neq\beta_3$, so away from $m=0$, which corresponds to helices; and we want a real square root, so we clip the radicand to be positive.  Then the least-squares solution is 
\begin{align}
    s_0=\frac{1}{\ell}\int_0^\ell \left(\frac{1}{\delta}F\left(\sqrt{\operatorname{max}\left(0,\frac{ \Big\langle\frac{-\lambda}{||\lambda||},\gamma'(Ss)\Big\rangle-\beta_3}{\beta_2-\beta_3}\right)}\Big|m\right)-s\right)\,\mathrm{d}s ,
\end{align}
with scale-invariant residual:
\begin{align}  \frac{1}{\ell} \int_0^\ell \left( \frac{1}{\delta}F\left(\sqrt{\operatorname{max}\left(0,\frac{ \Big\langle\frac{-\lambda}{||\lambda||},\gamma'(Ss)\Big\rangle-\beta_3}{\beta_2-\beta_3}\right)}\Big|m\right)-s-s_0\right)^2\,\mathrm{d}s.
\end{align}
Near the turning points of the pendulum, so the extrema of \ref{height}, the inverse \ref{inverse} is ill-conditioned. So this linear last-squares solution is only stable away from these extrema. \\

The choice of representative for $s_0$ modulo the period $2K(m)/\delta$ also affects the rotation and translation in the next step.  Different choices may 
be considered depending on the application.  The resulting elastic curve is always the same regardless of the choice, but  one should avoid a jump in the parameters if a continuous family of curves is to be constructed via interpolation of parameters.

\subsection{Recovery of rotation and translation}\hfill\\
Lastly, we need to determine the rotation $R$ and translation $\textbf{x}_0=(x_0,y_0,z_0)$. This can be done via the continuous version of the classical Kabsch algorithm. To summarize it:\\

We want to align the curve $\gamma(s)$ with $\zeta(s)=S\zeta_{\hat{\lambda}_0,\hat{\omega}}\left(\frac{s}{S}+s_0\right)$ for $s \in [0,L]$. Let $\bar{\gamma}$ and $\bar{\zeta}$ be the centroids of the two curves. We will translate both centroids to the origin and rotate. Thus, the translation $\textbf{x}_0=\bar\gamma-R\bar{\zeta}$ where $R$ is the rotation obtained from the following minimization problem:
\begin{align}
    \underset{R}{\text{minimize}}\int_0^{L}||(\gamma(s)-\bar{\gamma})-R(\zeta(s)-\bar{\zeta})||^2\,\mathrm{d}s. \end{align}
In fact, the rotation can also be determined slightly more efficiently, by first aligning the vector $(0,0,-1)$ corresponding to $\zeta_{\hat{\lambda}_0,\hat{\omega}}$ with the $\lambda$ recovered for $\gamma$, and then solving a one-dimensional minimization problem to find the rotation angle about the direction of $\lambda$.

%%%%%%%%%%%%%%%%%%%%%%%%%%%%%%%%%%%%%%%%%%%%%%%%%%%%%%%%%%%%%%%%%%

\section{Interactive design of elastic curve configurations using the least-squares approximation} 
\label{Design}

The least-squares parameter recovery method described in Section \ref{lssection} can be adapted to various applications where exact elastic curves that adhere approximately to a given shape are required to be produced at interactive speeds.  We consider some examples in this section.

\subsection{Interactive design of elastica-swept surfaces}
Design of elastica-swept surfaces using \emph{planar} elastica has been 
studied in \cite{brander2016designing, brander2018designing, fisker2019surface}, motivated by robotic hot-blade cutting (see Section \ref{hotbladecutting} below).
Here we consider how to extend this to
3D elastic curves.

One simple way to design elastica-swept surface patches is to take an arbitrary parameterized surface $S(u,v)$ as input, choose a finite subset of either the $u$ or $v$ parameter lines as a family of sweeping curves, then use the least-squares method to approximate each of these parameter lines by an elastic curve.  The fidelity to the input surface can be measured by, for instance, the maximum value of the $L^2$ difference between the target curve and the approximating curve among the chosen curves. Of course, other metrics such as an $H^1$ or $H^2$ metric could also be of interest here.

\begin{figure}[ht] 
\centering
\begin{small}
 \begin{tabular}{cc}
         \begin{tabular}{cc}
            \includegraphics[height=0.31\textheight]{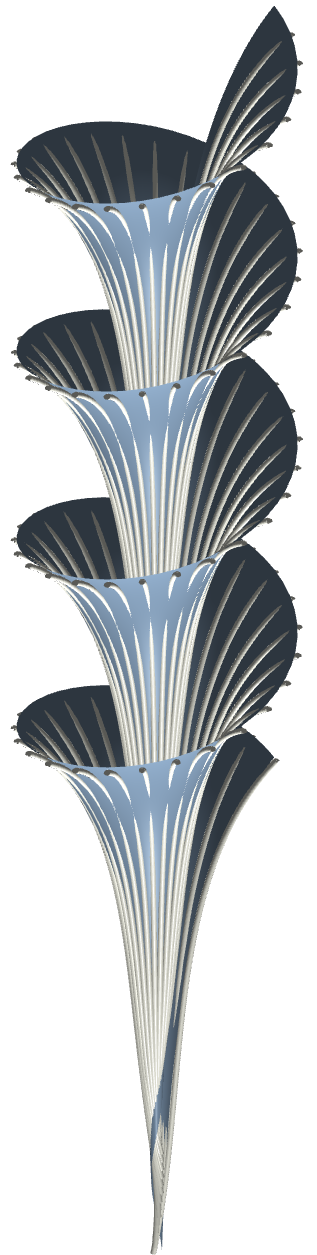} \,
             &  \,
        \includegraphics[height=0.31\textheight]{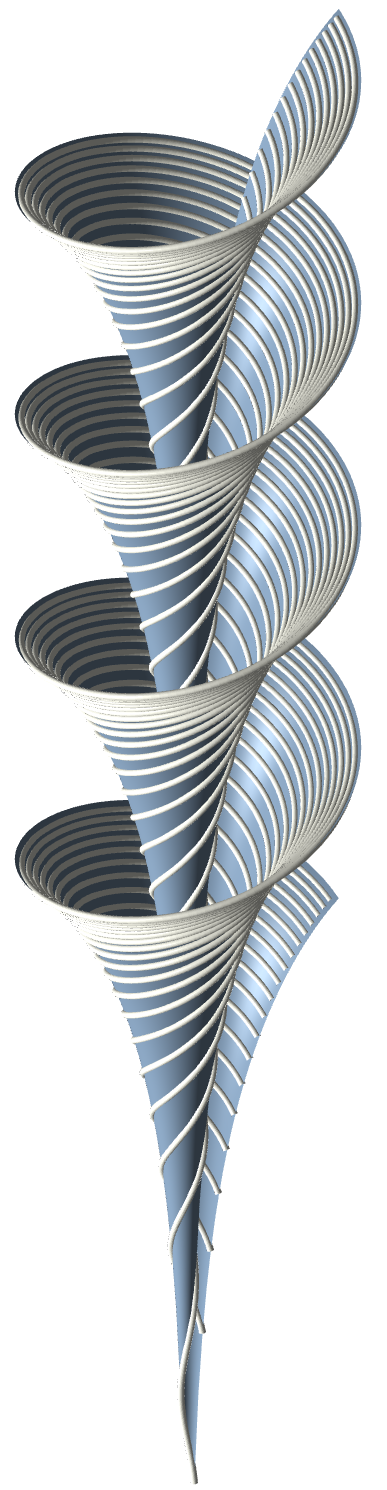} \vspace{1ex} \\
          Dini, Err $0.0056$ &   Dini, Err $0.0006$ 
            \end{tabular} \quad
     &     
           \begin{tabular}{cc}
                 \includegraphics[height=0.12\textheight]{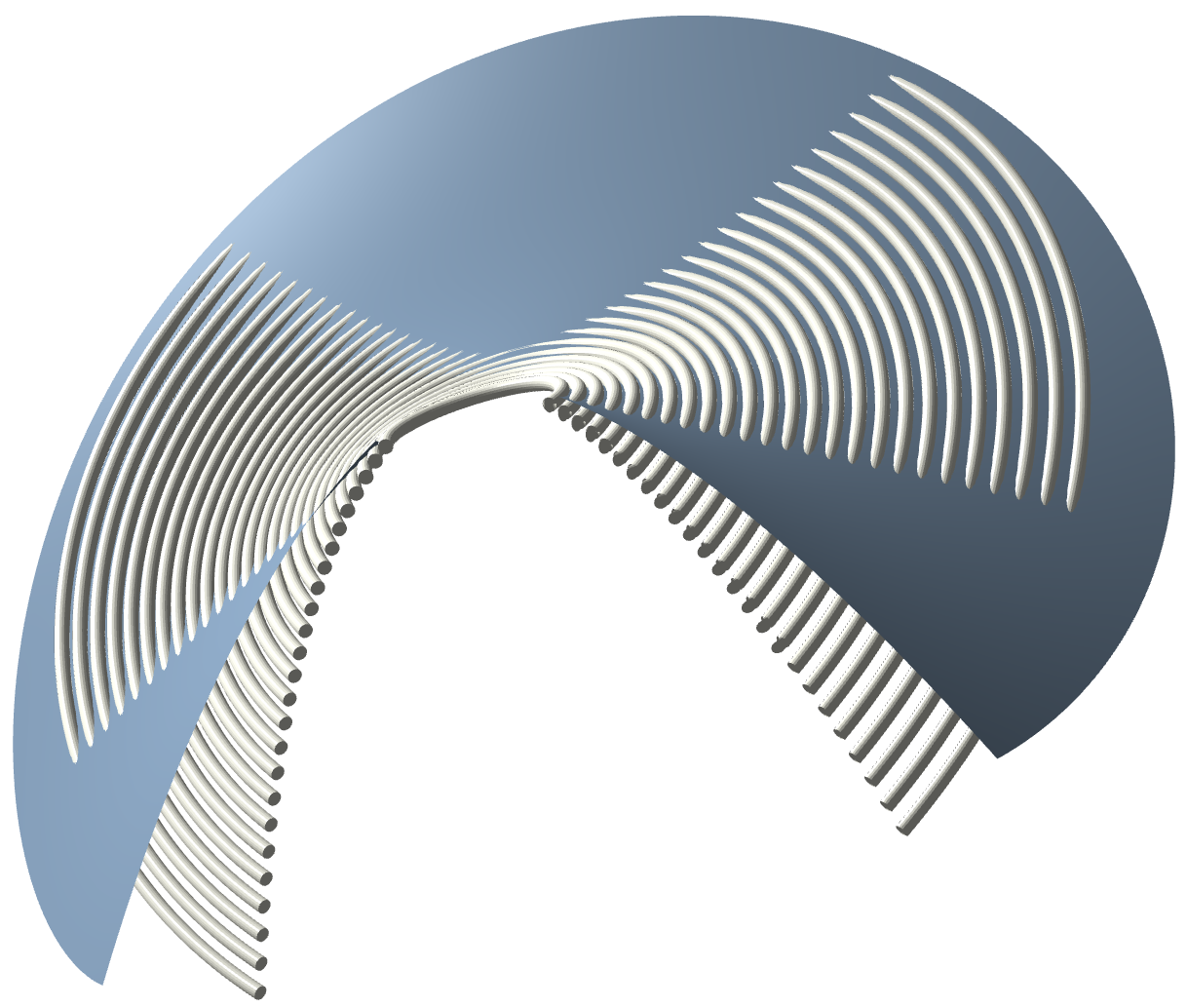} 
                 &   
     \includegraphics[height=0.12\textheight]{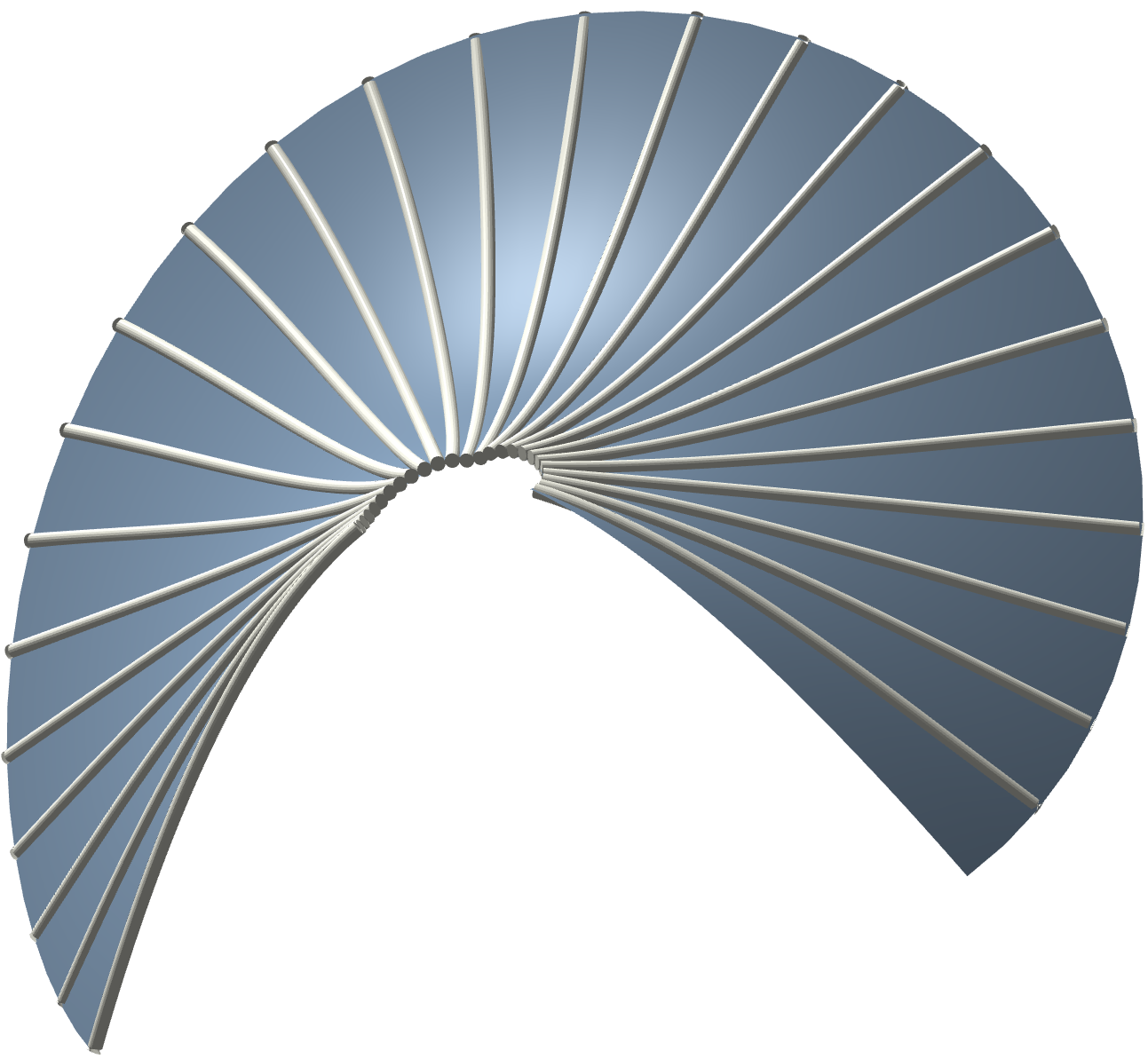} \vspace{2ex} \\
               Bour, Err $0.023$ & Bour, Err $0.0053$ \vspace{3ex} \\
      \includegraphics[height=0.1\textheight]{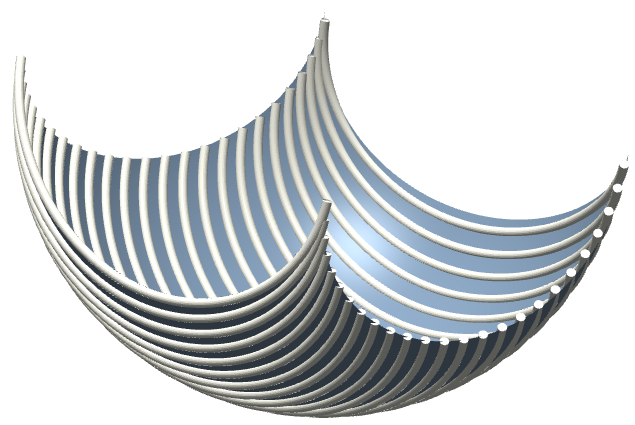}  &   
        \includegraphics[height=0.12\textheight]{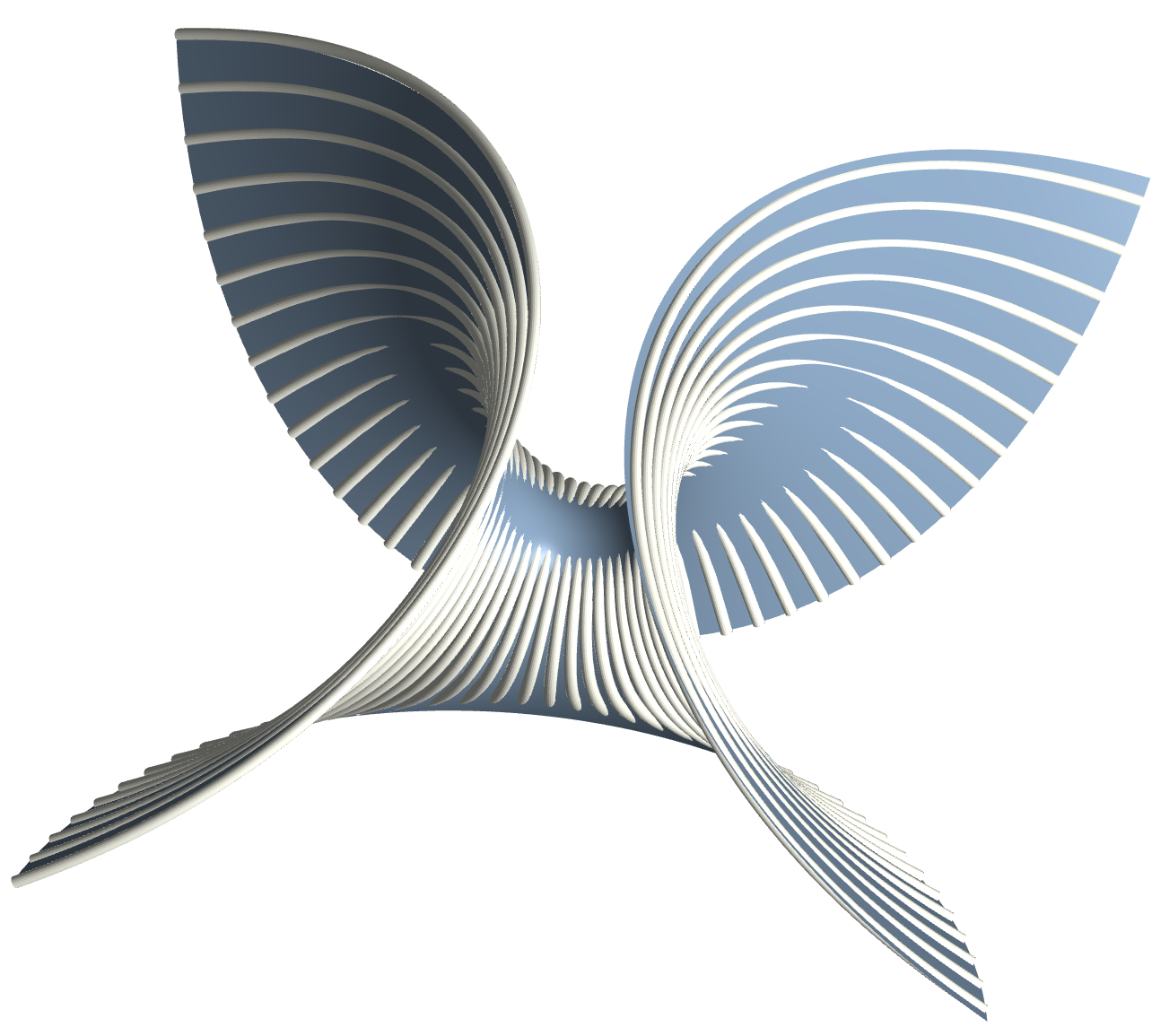}  \vspace{2ex} \\
            Sphere, Err $0.0045$ &   Enneper, Err $0.01$
            \end{tabular}  
    \end{tabular}  \vspace{1ex}\\
    \end{small}
  \caption{Examples of surfaces approximated by elastica-swept surfaces.}
      \label{figure2} 
\end{figure}

  Some examples of elastica-swept surfaces produced
  using this method are shown in
Figure  \ref{figure2}. The blue surfaces are the input surfaces and the white curves are plots of the least-squares elastic curves obtained from one of the parameter line families.  The error shown on the examples is the maximum among the curves of the (scaled) $L^2$ distance between the input curves $\gamma$ and the corresponding approximating elastic curves $\Gamma(t;\bsigma_\gamma)$, given
 by the scale-invariant quantity
\beq \label{err2}
\Err_2= \left(L^{-2}\int_0^1 \left|\left| \gamma(t) -\Gamma\left(\frac{s(t)}{L};\bsigma_\gamma\right)\right|\right|^2 \frac{s'(t)}{L} \textup{d}t\right)^{1/2},
\eeq
where $L$ is the length of the input curve $\gamma$ and $s(t)$ its arclength function.  Recall that $\Gamma(t;\bsigma_\gamma)$ has a constant speed parametrization. We want to parameterize it in such a way that it has the same speed as $\gamma$, for a meaningful comparison of the two curves.

For both Dini's pseudospherical surface and Bour's minimal surface, the same target surface patch is approximated using two different families of isolines: obviously the choice of parameter lines for an input surface can give different results for the output surface.  A sphere can, of course, be parameterized by families of great circles, which are elastica, but we have used a different parameterization by non-planar curves in our example.

\subsection{Robotic control of a linear element}
The problem of robotically guiding a flexible linear element through an obstacle
course has been studied in, for instance,  \cite{moll2006path,levin2025dual}.
A simple and fast procedure for determining a 
solution to this problem is as follows:
\begin{itemize}
\item \textbf{Step 1:} use any method to find a smooth simply connected surface that avoids all obstacles, with a buffer distance to allow for inaccuracies in the least-squares approximation.
\item \textbf{Step 2:} Approximate the surface by an elastica-swept surface. 
\end{itemize}
The above can be done at interactive speeds, so it is easy to verify the result and change the buffer sizes if the solution is not satisfactory.
Since the solution is given in terms of exact elastic curves, if a rod of length $L$ is to be guided through the course with robotic controllers placed at the endpoints, then the positions and tangents for the controllers can be read directly off the parameterized  curves. 

%\begin{figure*}[h]  
%\centering   \includegraphics[height=0.15\textheight]{images/obstacles1a} \quad
%\includegraphics[height=0.15\textheight]{images/obstacles1b}  \quad
%\includegraphics[height=0.15\textheight]%{images/obstacles1c} 
  %  \caption{Pathfinding for robotic steering of flexible linear element.}
%    \label{pathfinding}
%\end{figure*}

%%%%%%%%%%%%%%%%%%%%%%%%%%%%%%%%%%%
\section{The first residual $\Res_\lambda$ as a measure of input curve quality}\label{Residual}
For many applications (see, for example, the discussion of splines in Section \ref{splinesection} below) it is desirable to estimate how close the input curve $\gamma$ and the least-squares approximating elastica $\Gamma(t;\bsigma_{\gamma})$ are, without having to recover all the parameters and evaluating the difference between the curves.  

Regarding the input curve as simply a guide for the general shape of the desired elastic curve, we can, as mentioned above, first replace the input by a cubic spline.  As a general principle we expect that the geometric correlation between the least-squares approximation and the cubic spline, as measured by the first residual, can be well understood by studying a representative sample of cubic B\'ezier curves. This is a natural test space: cubic B\'ezier curves are common in geometric design, and their parameter space has dimension $12$, only one higher than that of the space of elastic curve segments. %reference to Ann-Sofie?

\begin{figure} [htb]
	\begin{center}	
\includegraphics[height=0.13\textheight]{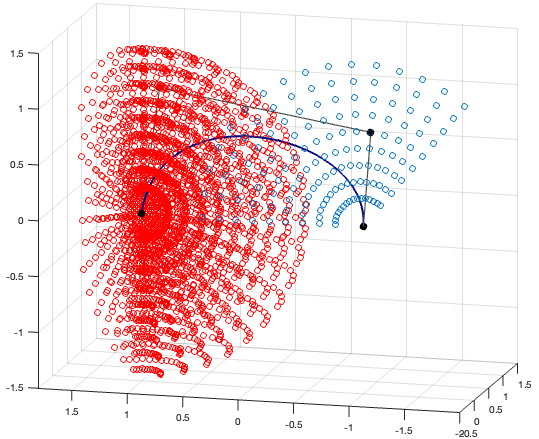}  	\quad \quad
\includegraphics[height=0.13\textheight]{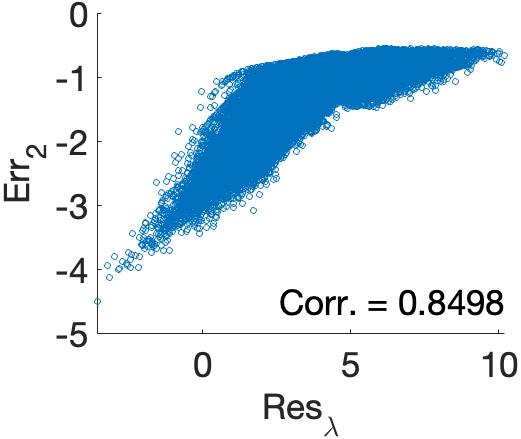}  
		\end{center}
\caption{Left: the sample set of middle control points. Right: log plot of
the $L^2$ error between the target curve and the elastica versus $\Res_\lambda$.} 
\label{sampleset}
\end{figure}

We  applied the first guess algorithm to a set of $435,157$ cubic B\'ezier curves, chosen to cover the range of shapes that are expected to arise as a segment of a spline in geometric design. The control points for the B\'ezier curves have the fixed outer control points $p_0=[0,-1,0]$,  $p_3 = [0, 1,0]$, and the two middle points $p_1$ and $p_2$ are given by all possible
pairs $(p_1,p_2)$ where $p_1$ is one of the blue points and $p_2$ one of the red,  in the sample set shown in Figure \ref{sampleset} (left). 
\begin{align*}
p_1 =& (0,y,z)=(0,-1+r \cos \phi, r \sin \phi), 
 \quad \quad (r,\phi) \in [0.25,1.5] \times  \left[0, \frac{3\pi}{4} \right], \\
p_2 =& (x,y,z)=(r\cos \theta \cos \phi, 1 + r \sin \theta \cos \phi, r \sin \phi),  \\
&
   \quad \quad (r, \theta, \phi) \in [0.24,1.5] \times \left[\frac{\pi}{2},\frac{5 \pi}{4}\right] \times \left[\frac{\pi}{2}+0.26, \frac{3\pi}{2}\right],
\end{align*}
and the parameter values are chosen uniformly from those intervals.  Any point pairs where the two points are within a distance of $0.2$ of each other were removed from the set, as these are not generally close to elastic curves.

The scatter plot in Figure \ref{sampleset} (right) shows the correlation between $L^2$ error $\Err_2$ of the resulting least-squares elastica and the target curve and the first residual. As the plot shows, they are well correlated with a Pearson correlation coefficient of $0.85$.

Since the first least-squares residual can be computed much more quickly than going through all of the steps in finding the least-squares elastica and then comparing that with the original curve, we would like to use this as a measure for how close
a potential target curve is to an elastic curve.
$10\%$ of curves have  $\Res_\lambda$ less than $75.8$, and of these the maximum $L^2$ error is $0.020$.
Likewise, $5\%$ have $\Res_\lambda$  less than $28.7$, and of these the maximum $L^2$ error is $0.016$.

\begin{figure}[htb]
\unitlength = 0.135\textheight
	\begin{center}	
    \includegraphics[height=\unitlength]{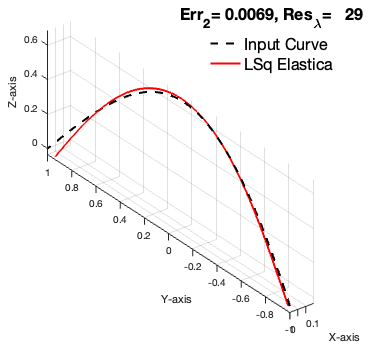}  \,
    \includegraphics[height=\unitlength]{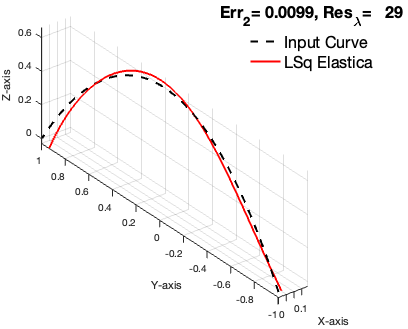} \, 
\includegraphics[height=\unitlength]{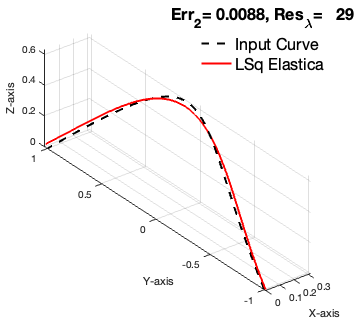}	\,  
\includegraphics[height=\unitlength]{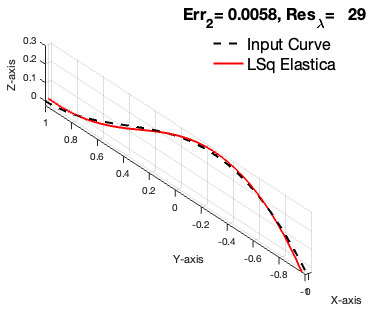}	 \vspace{2ex}\\
%%%%%%%%%%%%%%%%%%%%%%%
%\includegraphics[height=\unitlength]{images/10pc1.png}  \,
 %   \includegraphics[height=\unitlength]{images/10pc3.png} \, 
%\includegraphics[height=0.9\unitlength]{images/10pc4.png}	\,  
%\includegraphics[height=0.9\unitlength]{images/10pc7.png}	 \vspace{2ex}\\
%%%%%%%%%%%%%%%%%%%%%%%%%%%%%%%%%%%%%%%%
\includegraphics[height=\unitlength]{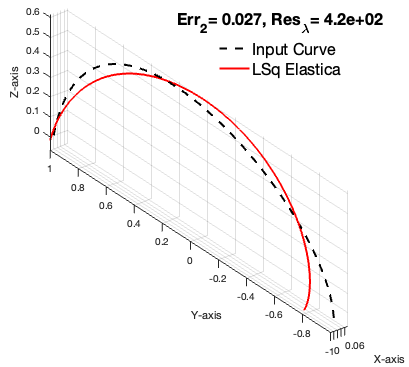}  \,
    \includegraphics[height=\unitlength]{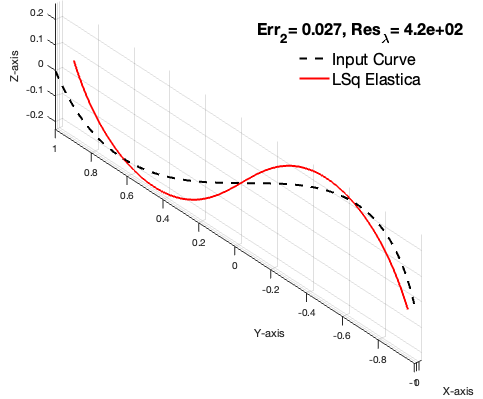} \, 
\includegraphics[height=\unitlength]{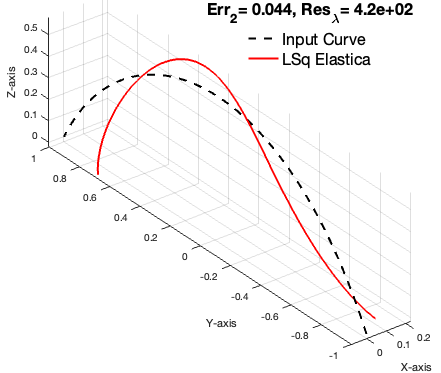}	\,  
\includegraphics[height=\unitlength]{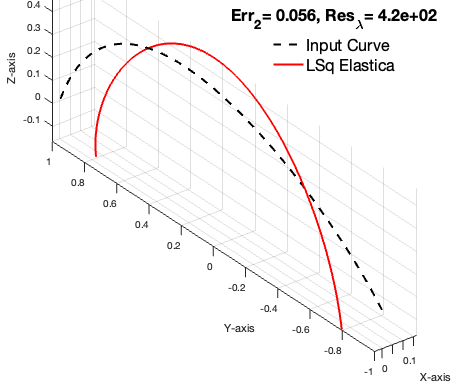}	 \vspace{2ex}\\
		\end{center}
\caption{Curves at the $5$th (top) and $25$th (bottom) percentiles of the least-squares residual.  }  
\label{fivepercentcurves}
\end{figure} 

 Curves at the $5$th and $25$th percentiles for $\Res_\lambda$ are shown in Figure \ref{fivepercentcurves}.
B\'ezier curves with similar values for $\Res_\lambda$ can be expected to be about as 
close to their least-squares approximating elastic curve as those depicted in this figure.

%%%%%%%%%%%%%%%%%%%%%%%%%%%%%%%%%%%%%%%%%%%%%%%%%%%
\section{Elastic Splines}  \label{splinesection}
For some applications an elastic spline is desirable. For example, cables can be routed with a series of connectors where the cable is free to rotate and/or slide through the joint, to create either $C^0$ or $C^1$ splines. Likewise, $C^0$ splines appear in fabrication of objects from thin flexible linear strips.  We consider here how this can be done using the least-squares elastica method.

\subsection{$C^0$ continuous splines}  
Converting an input curve into an elastic spline with $C^0$ continuity is simple using our method. It is a matter of dividing the curve into $n$ segments, then, for each segment, find the approximating elastica and finally scale, rotate and translate it so that the endpoints match those of the target curve segment. There is a freedom of rotation about the axis between the endpoints, and this rotation can be chosen by, for example, requiring that the rescaled curve is as close as possible to the input curve.

\begin{figure*} [hb]
	\begin{center}	
\includegraphics[height=0.20\textheight]{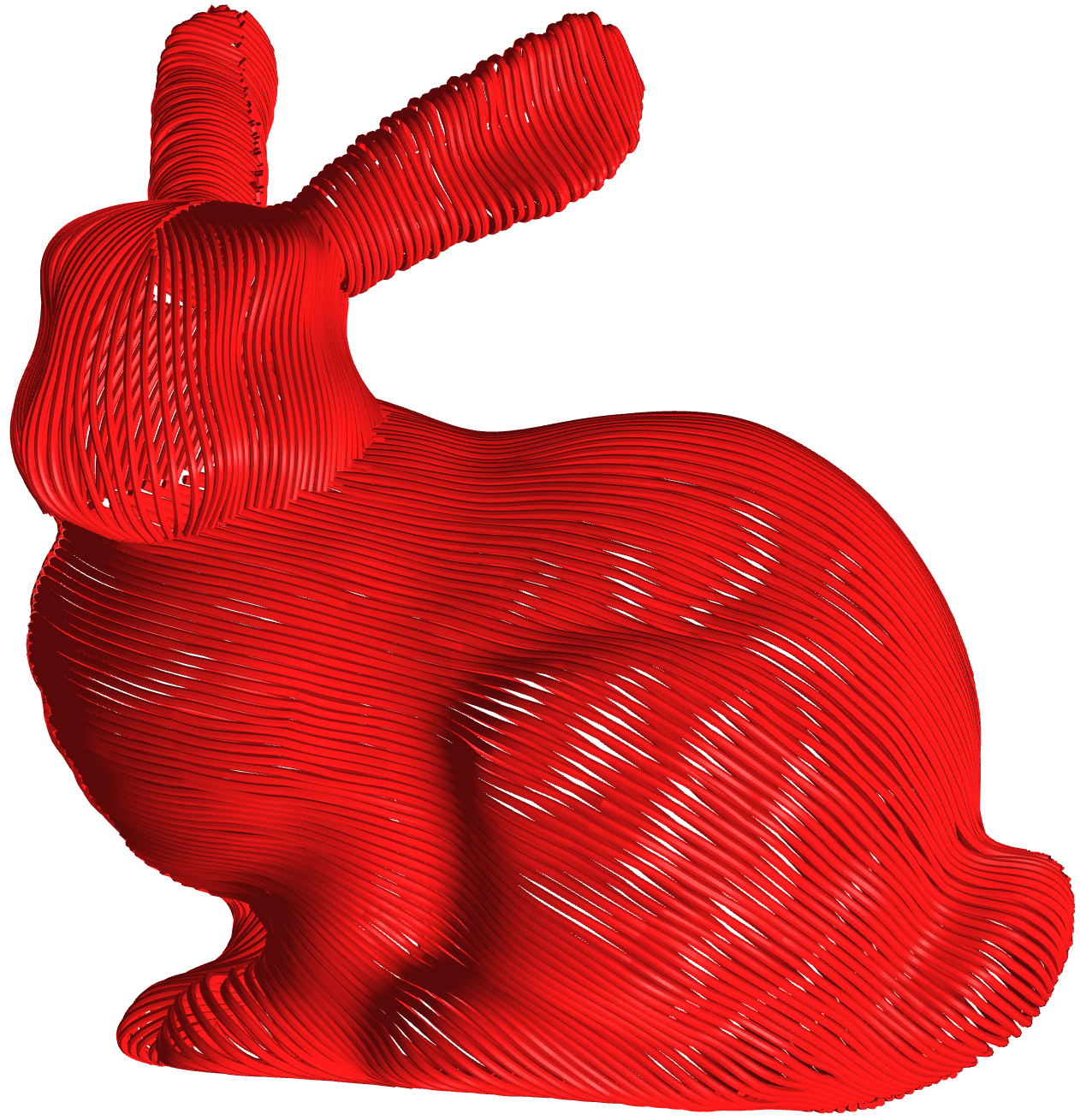}  	\quad\quad
\includegraphics[height=0.20\textheight]{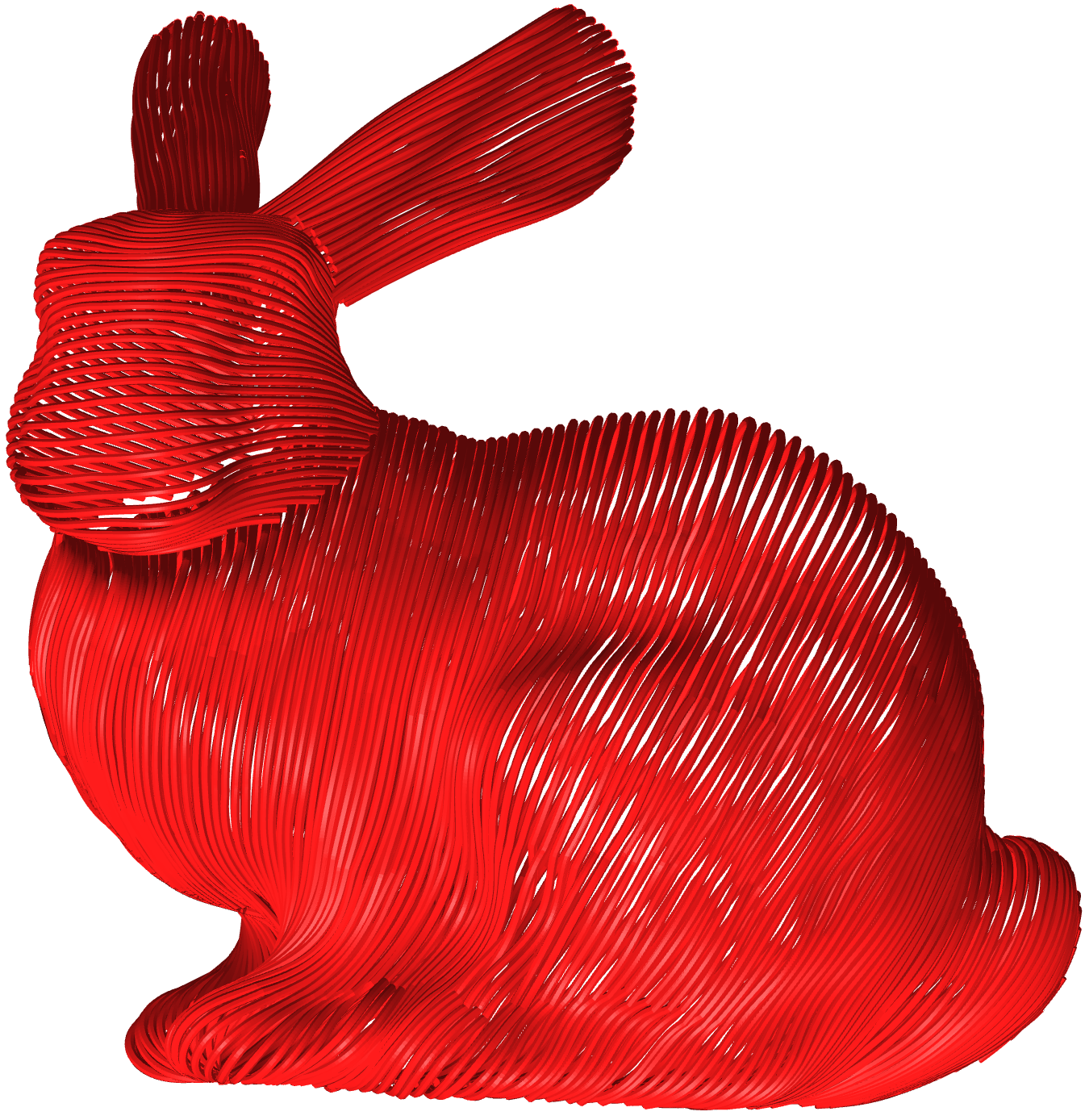}  
		\end{center}
\caption{Rationalizations of the bunny with $C^0$ elastic splines.} 
\label{bunny}
\end{figure*} 

An example is shown in Figure \ref{bunny}.  We started with a mesh of the Stanford bunny composed of four UV-patches\footnote{Obtained from {\fontfamily{cmss}\selectfont https://blenderartists.org/t/uv-unwrapped-stanford-bunny-happy-spring-equinox/1101297 }}. We extracted $u$- and $v$-polylines, and approximated those by smoothing splines, and in turn replaced each of those with an elastic $C^0$ spline.

Note that, if we want to find an elastic spline that is very close to the target curve, we can compute the first residual $\Res_\lambda$ iteratively for a segment of increasing length of the curve until it reaches a desired upper limit, thus producing a spline with $n$ segments with an upper bound on $\Res_\lambda$.

\subsection{$C^1$ continuous splines} 
\label{c1section}
A generalization can also be used to produce tangent continuous splines as follows: if two segments are required, first divide the curve into two segments. Find the least-squares elastica for each of segments $1$ and $2$. Attach segment $2$ to segment $1$ aligning their tangent vectors. Again there is a freedom of rotation about the tangent, and this rotation can be chosen by, for example, requiring that the Frenet frames match at the join, or such that the curve is as close as possible to the input curve. If interpolation of the endpoints is required, then this curve is again scaled, rotated and translated to match the original curve's endpoints.  The same procedure can be applied to any number of segments. Obviously, we can expect the final elastic spline will deviate further from the input curve the more translations and rotations are used.

\subsection{Using the least-squares residual to find intelligent input}
As has been done with planar elastic curves \cite{brander2018bezier,brander2018designing}, 
we can improve the utility of the methods mentioned above by working with intelligently chosen input.  For example, if the input curve is already close to an elastic curve, (or an elastic spline), then a $C^0$ elastic spline constructed as above can be expected to have essentially the same level of tangent continuity as the input curve.   Similarly, if we use the $C^1$ method described in Section
\ref{c1section}, we can expect the result to be close to the input curve.
\begin{figure}[hbt!]
\centering
\includegraphics[height=0.24\textheight]{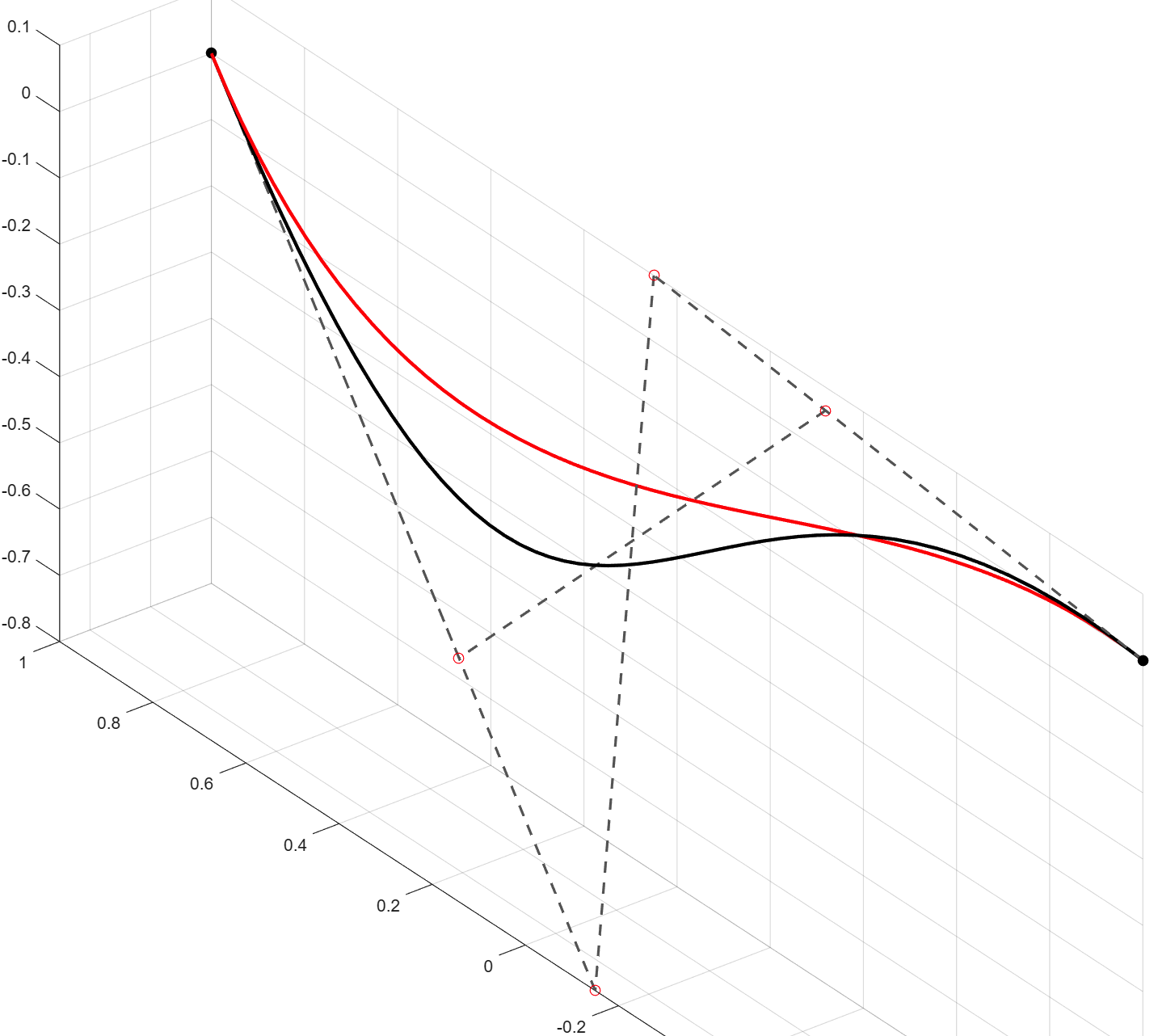}
\caption{Projecting a cubic  B\'ezier curve (black) to a cubic B\'ezier curve (red) that is closer to an elastic curve by moving the inner control points along the tangent directions to find a minimum $\Res_\lambda$ of $4.5$. }
\end{figure}

One can utilize the least-squares residual to work with B\'ezier curves and adjust the control points on the fly, as in \cite{brander2018bezier}.
A cubic B\'ezier curve is determined by prescribing a \textit{control polygon}, which is composed of the endpoints and two inner control points. Moving the two inner control points up and down the tangent directions will change the shape without changing the tangent at the endpoints. We can compute automatically, at interactive speeds, the optimum amount to adjust these to minimize $\Res_\lambda$. Applying the same procedure to a cubic spline will convert the cubic spline into a spline where each segment has a reduced $\Res_\lambda$.  Depending on the application, we can set the maximum allowable value for $\Res_\lambda$ to decide whether or not our spline is acceptable as input for our methods.
In fact, if we modify an interpolating cubic spline, the new curve still interpolates the desired points and passes through them at the same tangent directions. The elastic spline obtained from replacing a cubic spline with low $\Res_\lambda$ can then be interactively controlled by moving the same control points of the original spline.

%%%%%%%%%%%%%%%%%%%%%%%%%%%%%%%%%%%%%%%%%%%%%%%%%%%%%%%%%%
\section{Approximation of a given target curve by an elastic curve}\label{Approximation}
\subsection{Nonlinear optimization}\hfill\\
The first guess can be used as an initial value for a range of optimization problems. Here we are trying to approximate a given unit-speed curve $\gamma$ by an elastic segment of the same length, so the first guess $\Gamma(t;\bsigma_{\gamma})$ is refined by minimizing the scale-invariant error
\begin{align*}
  \left(\frac1{L^3}
  \int_0^L
  \|\gamma(s)-\Gamma\left(\frac{s}{L};\bsigma\right)\|^2\,ds \right) ^{1/2}
\end{align*}
where
\[
  \bsigma=(\hat{\lambda}_0,\hat{\omega},s_0,\ell,S,R,\textbf{x}_0)
\]
and $\Gamma(t;\bsigma)$ is the corresponding elastic segment.

In the experiments, this minimization was performed using MATLAB's
\texttt{fmincon} with the SQP algorithm. The initial value was the least-squares elastica $\Gamma(t;\bsigma_{\gamma})$, and the admissibility of $(\hat{\lambda}_0,\hat{\omega})$ was enforced during the optimization.

\subsection{Numerical examples}
\begin{figure}[htb]
	\begin{center}	
        \includegraphics[height=33mm]{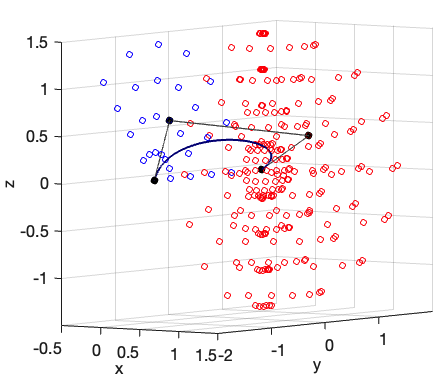}  \quad 
          \includegraphics[height=33mm]{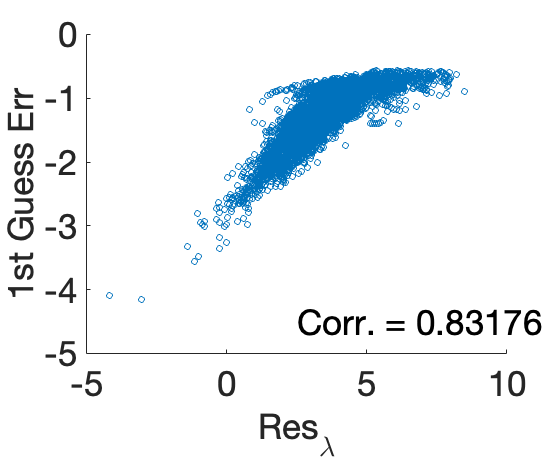} \,
           \includegraphics[height=33mm]{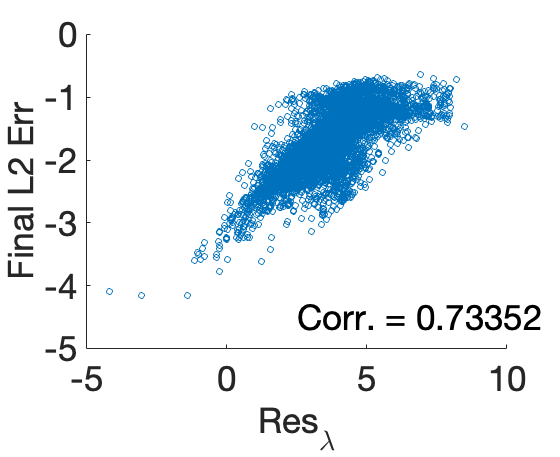}	
		\end{center}
\caption{Left: middle control point sample space. Right: log plots of least-squares residual versus error in approximation by first guess and final optimized solution.}  \label{figscatter}
\end{figure}

\hfill\\
We took a sample of 7035 cubic B\'ezier curves with outer control points $p_0 = (0,-1,0)$, $p_3=(0,1,0)$ and the two middle control points $p_1$ and $p_2$ chosen from (Figure \ref{figscatter}, left) all possible pairs:
\begin{align*}
p_1 &=& (0,y,z)=(0,-1+r \cos \phi, r \sin \phi), 
 \quad \quad (r,\phi) \in [0.3,1.5] \times  [0, 3\pi/4], \\
p_2 &=& (x,y,z)=(r\cos \theta \cos \phi, 1 + r \sin \theta \cos \phi, r \sin \phi),  \\
&&
   \quad \quad 
   (r, \theta, \phi) \in [0.3,1.5] \times [\pi/2, 5 \pi/4 ] \times  [\pi/2+0.26, 3\pi/2],
\end{align*}
and the values are chosen in intervals of $0.38$.  Curves where $|p_2-p_3|<0.3$ were removed from the sample.

The B\'ezier curves in the sample were approximated using  the first guess followed by an optimization using MATLAB's fmincon with the SQP algorithm as described above.  In Figure \ref{figscatter}
we have plotted log plots of the error for both the first guess and the optimized solution against the least-squares residual $\operatorname{Res}_{\lambda}$ given by \eqref{Eres}.   \\

Note that the residual $\operatorname{Res}_{\lambda}$, which is fast to compute, is a good indicator for how well the input curve can be approximated by an elastic curve, with correlation coefficients of $0.83$ and $0.73$ for the errors in the first guess and optimized solution respectively.
\begin{figure}[htb]
	\begin{center}
    $
    \begin{array}{ccc}    \includegraphics[height=33mm]{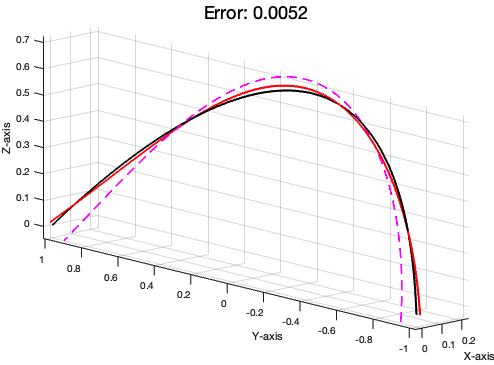}   &
    \includegraphics[height=33mm]{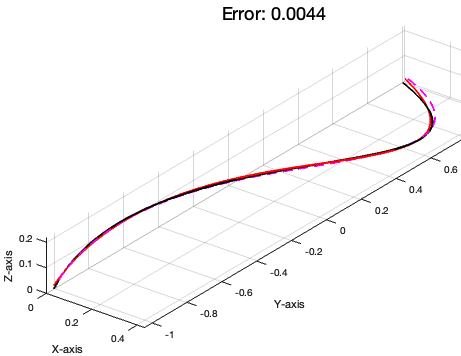} 
	&
  \includegraphics[height=33mm]{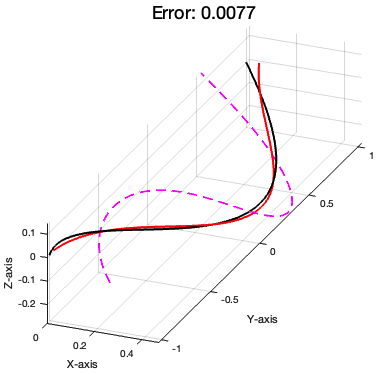}  \vspace{3ex}\\   
    \includegraphics[height=33mm]{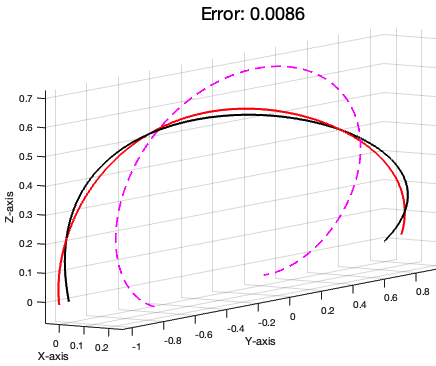} &
      \includegraphics[height=33mm]{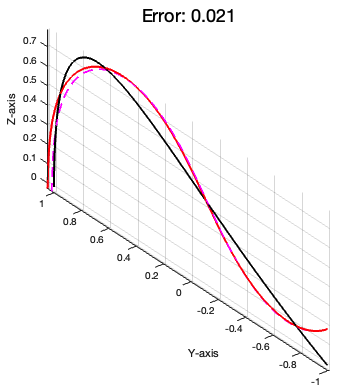}  & \quad
    \includegraphics[height=33mm]{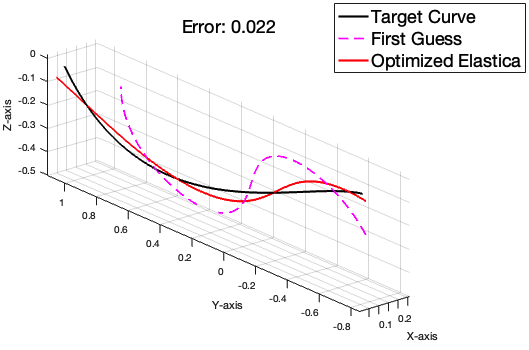}
    \end{array}
$
           \end{center}
\caption{Approximations of cubic B\'ezier curves by 3D elastic curves, with approximation error at the 10th percentile (left two), 25th percentile (middle two), and 50th percentile (right two) in the sample space.}  \label{approximated}
\end{figure} 
The mean and median for the first guess error were $0.092$ and $0.0784$ and these improved to, respectively, $0.035$ and $0.021$ for the optimized solution.  The top quartile of the optimized curves have an error of $0.0086$ or less, and the top $10\%$ have an error at most $0.0044$.  Examples of these are shown in Figure \ref{approximated}.

\subsection{Application: Surface rationalization for hot-blade cutting}
\label{hotbladecutting}
Surface rationalization for hot-blade cutting was studied in 
\cite{norbjerg2017rationalization, sondergaard2016robotic,petkov2017hot} using planar elastic curves.
Using 3D elastic curves has the advantage of further degrees of freedom, allowing a closer approximation to a given target surface.

The general problem involves first segmenting the target surface into patches that are small enough to fit inside a block of EPS and to be well approximated by a family of elastic curve segments. Then a suitable parameterization should be given on each patch such that one family of parameter lines is well approximated by elastic curve segments. 

The global problem is the subject of ongoing work. The local problem of approximating a single patch is well addressed by the method described in this paper.
An example is shown in Figure \ref{FigSurfApprox} to the right.

\section{Conclusion and Future Directions} \label{conclusion}
We have presented a stable parameter recovery method for three-dimensional elastic curve segments and used it to approximate general space curves. The method can be used either as a fast way to produce elastic curves for various applications, or as a first guess for nonlinear optimization. Whereas, in this paper, we have optimized for closeness in the $L^2$ sense and it is straightforward to use the same method with, for example, an $H^1$-metric,  we expect to be able to use the same first guess for other optimization problems, for instance solving boundary-value problems.

The least-squares elastica method is a very fast way to produce exact 3D elastic curves guided by the \emph{shape} of a target curve. With an implementation in MATLAB that has not been optimized for speed, the surfaces in Figure \ref{figure2}, typically consisting of about $50$-$100$ elastic curves are computed in a few hundredths of a second. The Stanford bunny rationalization, consisting of around 6000 elastic curves, is computed in about 5 seconds on an ordinary laptop.
Design tools can thus be made to work with very large numbers of elastic curves interactively. 

Applications include geometric modeling with physically realizable curves, such as flexible strips, cable routing, bent rods and hot-blade cutting.  For elastic rods that include a mechanical twist, we expect that a similar method can be developed where the extra relevant parameter, the scaled physical twist density (see the formulation in \cite{ameline2017classifications}), should also be obtainable by a least-squares method.
Other directions for future work include
piecewise-elastica approximation of long curves, and global approximation of surfaces by multi-patch elastica-swept surfaces.

\section*{Code}
MATLAB implementations of the least-squares approximation as well as the optimized approximation of a curve by an elastica can be found at this
repository: 
\[
\textup{\href{https://github.com/mishadtu/3d-elastic-curves}{https://github.com/mishadtu/3d-elastic-curves}}
\]

\section*{Acknowledgements}
Research supported by Horizon-EIC-2023 Pathfinder Challenge 01 project nr. \ 101161085
STACK.

\bibliographystyle{amsplain} 
 \bibliography{mybib}

@article{jin2022robotic,
  title={Robotic cable routing with spatial representation},
  author={Jin, Shiyu and Lian, Wenzhao and Wang, Changhao and Tomizuka, Masayoshi and Schaal, Stefan},
  journal={IEEE Robotics and Automation Letters},
  volume={7},
  number={2},
  pages={5687--5694},
  year={2022},
  publisher={IEEE}
}

@inproceedings{monguzzi2024potential,
  title={Potential field-based online path planning for robust cable routing},
  author={Monguzzi, Andrea and Mantegna, Niccol{\`o} and Zanchettin, Andrea Maria and Rocco, Paolo},
  booktitle={2024 IEEE/RSJ International Conference on Intelligent Robots and Systems (IROS)},
  pages={7558--7564},
  year={2024},
  organization={IEEE}
}

@article{levin2025dual,
  title={Dual Arm Steering of Flexible Linear Objects in 2-D and 3-D Environments Using Euler's Elastica Solutions},
  author={Levin, A and Grinberg, I and Rimon, ED and Shapiro, A},
  journal={IEEE Robotics and Automation Letters},
  year={2025},
  publisher={IEEE}
}

@article{moll2006path,
  title={Path planning for deformable linear objects},
  author={Moll, Mark and Kavraki, Lydia E},
  journal={IEEE Transactions on Robotics},
  volume={22},
  number={4},
  pages={625--636},
  year={2006},
  publisher={IEEE}
}

@incollection{sondergaard2016robotic,
  title={Robotic hot-blade cutting: an industrial approach to cost-effective production of double curved concrete structures},
  author={S{\o}ndergaard, Asbj{\o}rn and Feringa, Jelle and N{\o}rbjerg, Toke and Steenstrup, Kasper and Brander, David and Graversen, Jens and Markvorsen, Steen and B{\ae}rentzen, Andreas and Petkov, Kiril and Hattel, Jesper and others},
  booktitle={Robotic fabrication in architecture, art and design 2016},
  pages={150--164},
  year={2016},
  publisher={Springer}
}

@article{norbjerg2017rationalization,
  title={Rationalization in architecture with surfaces foliated by elastic curves},
  author={N{\o}rbjerg, Toke Bjerge},
  year={2017},
  publisher={Technical University of Denmark}
}

@article{petkov2017hot,
  title={Hot-blade cutting of EPS foam for double-curved surfaces—numerical simulation and experiments},
  author={Petkov, Kiril P and Hattel, Jesper H},
  journal={The International Journal of Advanced Manufacturing Technology},
  volume={93},
  number={9},
  pages={4253--4264},
  year={2017},
  publisher={Springer}
}

@phdthesis{fisker2019surface,
title = "Surface design and rationalization for robotic hot-blade cutting",
author = "Ann-Sofie Fisker",
year = "2019",
language = "English",
series = "DTU Compute PHD-2018",
publisher = "DTU Compute",
}

@article{brander2018designing,
  title={Designing interactively with elastic splines},
  author={Brander, David and B{\ae}rentzen, Jakob Andreas and Fisker, Ann-Sofie and Gravesen, Jens},
  journal={Computer Aided Geometric Design},
  volume={62},
  pages={181--191},
  year={2018},
  publisher={Elsevier}
}

@inproceedings{brander2016designing,
  title={Designing for hot-blade cutting: geometric approaches for high-speed manufacturing of doubly-curved architectural surfaces},
  author={Brander, David and B{\ae}rentzen, Jakob Andreas and Clausen, Kenn and Fisker, Ann-Sofie and Gravesen, Jens and Lund, Morten N and N{\o}rbjerg, Toke Bjerge and Steenstrup, Kasper Hornbak and S{\o}ndergaard, Asbj{\o}rn},
  booktitle={Advances in architectural geometry (AAG 2016)},
  pages={306--327},
  year={2016},
  organization={vdf Hochschulverlag AG an der ETH Z{\"u}rich}
}

@article{brander2018bezier,
  title={B{\'e}zier curves that are close to elastica},
  author={Brander, David and B{\ae}rentzen, Jakob Andreas and Fisker, Ann-Sofie and Gravesen, Jens},
  journal={Computer-Aided Design},
  volume={104},
  pages={36--44},
  year={2018},
  publisher={Elsevier}
}

@article{bgn2016,
  author  = {D. Brander and J. Gravesen and T. N\o{}rbjerg},
  title   = {Approximation by planar elastic curves},
  year    = {2016},
  journal = {Adv. Comput. Math. (2016)},
	volume={43},
	pages={25-43},
  doi     = {10.1007/s10444-016-9474-z}
}

@article{ameline2017classifications,
  title={Classifications of ideal 3D elastica shapes at equilibrium},
  author={Ameline, Olivier and Haliyo, Sinan and Huang, Xingxi and Cognet, Jean AH},
  journal={Journal of Mathematical Physics},
  volume={58},
  number={6},
  year={2017},
  publisher={AIP Publishing}
}

@book{pinkall2024differential,
  title={Differential Geometry: From Elastic Curves to Willmore Surfaces},
  author={Pinkall, Ulrich and Gross, Oliver},
  year={2024},
  publisher={Springer Nature}
}

@inproceedings{singer2008lectures,
  title={Lectures on elastic curves and rods},
  author={Singer, David A},
  booktitle={AIP Conference Proceedings},
  volume={1002},
  number={1},
  pages={3},
  year={2008},
  organization={Citeseer}
}

@inproceedings{mumford1994elastica,
  title={Elastica and computer vision},
  author={Mumford, David},
  booktitle={Algebraic Geometry and its Applications: Collections of Papers from Shreeram S. Abhyankar’s 60th Birthday Conference},
  pages={491--506},
  year={1994},
  organization={Springer}
}

@incollection{bergou2008discrete,
  title={Discrete elastic rods},
  author={Bergou, Mikl{\'o}s and Wardetzky, Max and Robinson, Stephen and Audoly, Basile and Grinspun, Eitan},
  booktitle={ACM Siggraph 2008 Papers},
  pages={1--12},
  year={2008}
}

@article{horn1983curve,
  title={The curve of least energy},
  author={Horn, Berthold KP},
  journal={ACM Transactions on Mathematical Software (TOMS)},
  volume={9},
  number={4},
  pages={441--460},
  year={1983},
  publisher={ACM New York, NY, USA}
}

@article{shen2003euler,
  title={Euler's elastica and curvature-based inpainting},
  author={Shen, Jianhong and Kang, Sung Ha and Chan, Tony F},
  journal={SIAM journal on Applied Mathematics},
  volume={63},
  number={2},
  pages={564--592},
  year={2003},
  publisher={SIAM}
}

@article{hafner2021design,
  title={The design space of plane elastic curves},
  author={Hafner, Christian and Bickel, Bernd},
  journal={ACM Transactions on Graphics (TOG)},
  volume={40},
  number={4},
  pages={1--20},
  year={2021},
  publisher={ACM New York, NY, USA}
}

@inproceedings{lee2020ruled,
  title={From ruled surfaces to elastica-ruled surfaces: new possibilities for creating architectural forms},
  author={Lee, Ting-Uei and Xie, Yi Min},
  booktitle={Proceedings of IASS Annual Symposia},
  volume={2020},
  number={30},
  pages={1--12},
  year={2020},
  organization={International Association for Shell and Spatial Structures (IASS)}
}

@article{bi2023design,
  title={Design and construction of kinetic structures based on elastica strips},
  author={Bi, Minghao and He, Yunzhen and Li, Zhi and Lee, Ting-Uei and Xie, Yi Min},
  journal={Automation in Construction},
  volume={146},
  pages={104659},
  year={2023},
  publisher={Elsevier}
}

@article{arroyo2020boundary,
  title={Boundary value problems for Euler-Bernoulli planar elastica. A solution construction procedure},
  author={Arroyo, Josu J and Garay, {\'O}scar J and P{\'a}mpano, {\'A}lvaro},
  journal={Journal of Elasticity},
  volume={139},
  number={2},
  pages={359--388},
  year={2020},
  publisher={Springer}
}

@article{wang2017real,
  title={Real time simulation of inextensible surgical thread using a Kirchhoff rod model with force output for haptic feedback applications},
  author={Wang, Zhujiang and Fratarcangeli, Marco and Ruimi, Annie and Srinivasa, AR},
  journal={International Journal of Solids and Structures},
  volume={113},
  pages={192--208},
  year={2017},
  publisher={Elsevier}
}

@article{panneerselvam2020constrained,
  title={A constrained spline dynamics (CSD) method for interactive simulation of elastic rods},
  author={Panneerselvam, Karthikeyan and Rahul and De, Suvranu},
  journal={Computational mechanics},
  volume={65},
  number={2},
  pages={269--291},
  year={2020},
  publisher={Springer}
}

@article{armanini2023soft,
  title={Soft robots modeling: A structured overview},
  author={Armanini, Costanza and Boyer, Fr{\'e}d{\'e}ric and Mathew, Anup Teejo and Duriez, Christian and Renda, Federico},
  journal={IEEE Transactions on Robotics},
  volume={39},
  number={3},
  pages={1728--1748},
  year={2023},
  publisher={IEEE}
}

@article{choi2024dismech,
  title={Dismech: A discrete differential geometry-based physical simulator for soft robots and structures},
  author={Choi, Andrew and Jing, Ran and Sabelhaus, Andrew P and Jawed, Mohammad Khalid},
  journal={IEEE Robotics and Automation Letters},
  volume={9},
  number={4},
  pages={3483--3490},
  year={2024},
  publisher={IEEE}
}

@article{cuvilliers2018comparison,
  title={A comparison of two algorithms for the simulation of bending-active structures},
  author={Cuvilliers, Pierre and Yang, Justina R and Coar, Lancelot and Mueller, Caitlin},
  journal={International Journal of Space Structures},
  volume={33},
  number={2},
  pages={73--85},
  year={2018},
  publisher={SAGE Publications Sage UK: London, England}
}

@article{lazaro2018mechanical,
  title={Mechanical models in computational form finding of bending-active structures},
  author={L{\'a}zaro, Carlos and Bessini, Juan and Monle{\'o}n, Salvador},
  journal={International Journal of Space Structures},
  volume={33},
  number={2},
  pages={86--97},
  year={2018},
  publisher={SAGE Publications Sage UK: London, England}
}

@article{alessi2024rod,
  title={Rod models in continuum and soft robot control: a review},
  author={Alessi, Carlo and Agabiti, Camilla and Caradonna, Daniele and Laschi, Cecilia and Renda, Federico and Falotico, Egidio},
  journal={arXiv preprint arXiv:2407.05886},
  year={2024}
}

@article{tola2022simulation,
  title={A Simulation of an Elastic Filament Using Kirchhoff Model},
  author={Tola, Saimir and Daci, Alfred and Zavalani, Gentian},
  journal={Mathematics and Statistics},
  volume={10},
  number={1},
  pages={25--34},
  year={2022}
}

@article{da2020simple,
  title={A simple finite element for the geometrically exact analysis of Bernoulli--Euler rods},
  author={da Costa e Silva, C{\'a}tia and Maassen, Sascha F and Pimenta, Paulo M and Schr{\"o}der, J{\"o}rg},
  journal={Computational Mechanics},
  volume={65},
  number={4},
  pages={905--923},
  year={2020},
  publisher={Springer}
}

@article{carlson1995numerical,
  title={Numerical computation of real or complex elliptic integrals},
  author={Carlson, Bille C},
  journal={Numerical Algorithms},
  volume={10},
  number={1},
  pages={13--26},
  year={1995},
  publisher={Springer}
}

@article{batista2019elfun18,
  title={Elfun18--A collection of MATLAB functions for the computation of elliptic integrals and Jacobian elliptic functions of real arguments},
  author={Batista, Milan},
  journal={SoftwareX},
  volume={10},
  pages={100245},
  year={2019},
  publisher={Elsevier}
}

\end{document}